# HOW MANY ZEROS OF A RANDOM POLYNOMIAL ARE REAL?

ALAN EDELMAN AND ERIC KOSTLAN

Abstract. We provide an elementary geometric derivation of the Kac integral formula for the expected number of real zeros of a random polynomial with independent standard normally distributed coefficients. We show that the expected number of real zeros is simply the length of the moment curve $(1, t, \ldots, t^n)$ projected onto the surface of the unit sphere, divided by $\pi$. The probability density of the real zeros is proportional to how fast this curve is traced out.

We then relax Kac's assumptions by considering a variety of random sums, series, and distributions, and we also illustrate such ideas as integral geometry and the Fubini-Study metric.

## Contents

Received by the editors November 22, 1993, and, in revised form, July 25, 1994.
1991 *Mathematics Subject Classification.* Primary 34F05; Secondary 30B20.
Supported by the Applied Mathematical Sciences subprogram of the Office of Energy Research, U.S. Department of Energy, under Contract DE-AC03-76SF00098.

## 1. Introduction

What is the expected number of real zeros $E_n$ of a random polynomial of degree $n$? If the coefficients are independent standard normals, we show that as $n \to \infty$,

$$E_n = \frac{2}{\pi}\log(n) \; + \; 0.6257358072... \; + \; \frac{2}{n\pi} \; + \; O(1/n^2) \; .$$

The $\frac{2}{\pi}\log n$ term was derived by Kac in 1943 [26], who produced an integral formula for the expected number of real zeros. Papers on zeros of random polynomials include [3], [16], [23], [34], [41, 42] and [36]. There is also the comprehensive book of Bharucha-Reid and Sambandham [2].

We will derive the Kac formula for the expected number of real zeros with an elementary geometric argument that is related to the Buffon needle problem. We present the argument in a manner such that precalculus level mathematics is sufficient for understanding (and enjoying) the introductory arguments, while elementary calculus and linear algebra are sufficient prerequisites for much of the paper. Nevertheless, we introduce connections with advanced areas of mathematics.

A seemingly small variation of our opening problem considers random $n$th degree polynomials with independent normally distributed coefficients, each with mean zero, but with the variance of the $i^{\text{th}}$ coefficient equal to $\binom{n}{i}$ (see [4], [31], [46]). This particular random polynomial is probably the more natural definition of a random polynomial. It has

$$E_n = \sqrt{n}$$

real zeros on average.

As indicated in our table of contents, these problems serve as the departure point for generalizations to systems of equations and the real or complex zeros of other collections of random functions. For example, we consider power series, Fourier series, sums of orthogonal polynomials, Dirichlet series, matrix polynomials, and systems of equations.

Section 2 begins with our elementary geometric derivation. Section 3 shows how a large class of random problems may be covered in this framework. In Section 4 we reveal what is going on mathematically. Section 5 studies arbitrary distributions but focuses on the non-central normal. Section 6 relates random polynomials to random matrices, while Section 7 extends our results to systems of equations. Complex roots, which are ignored in the rest of paper, are addressed in Section 8. We relate random polynomials to the Buffon needle problem in Section 9.

## 2. Random polynomials and elementary geometry

Section 2.1 is restricted to elementary geometry. Polynomials are never mentioned. The relationship is revealed in Section 2.2.

### 2.1. How fast do equators sweep out area?

We will denote (the surface of) the unit sphere centered at the origin in $\mathbb{R}^{n+1}$ by $S^n$. Our figures correspond to the case $n = 2$. Higher dimensions provide no further complications.

**Definition 2.1.** If $P \in S^n$ is any point, the *associated equator* $P_\perp$ is the set of points of $S^n$ on the plane perpendicular to the line from the origin to $P$.

This generalizes our familiar notion of the earth's equator, which is equal to $(\text{north pole})_\perp$ and also equal to $(\text{south pole})_\perp$. See Figure 1. Notice that $P_\perp$ is always a unit sphere ("great hypercircle") of dimension $n-1$.

Let $\gamma(t)$ be a (rectifiable) curve on the sphere $S^n$.

**Definition 2.2.** Let $\gamma_\perp$, the *equators of a curve*, be the set $\{P_\perp | P \in \gamma\}$.

Assume that $\gamma$ has a finite length $|\gamma|$. Let $|\gamma_\perp|$ to be the area "swept out" by $\gamma_\perp$ (we will provide a precise definition shortly). We wish to relate $|\gamma|$ to $|\gamma_\perp|$.

If the curve $\gamma$ is a small section of a great circle, then $\cup\gamma_\perp$ is a lune, the area bounded by two equators as illustrated in Figure 2. If $\gamma$ is an arc of length $\theta$, then our lune covers $\theta/\pi$ of the area of the sphere. The simplest case is $\theta = \pi$. We thus obtain the formula valid for arcs of great circles, namely,

$$\frac{|\gamma_\perp|}{\text{area of } S^n} = \frac{|\gamma|}{\pi}.$$

If $\gamma$ is not a section of a great circle, we may approximate it by a union of small great circular arcs, and the argument is still seen to apply.

The alert reader may notice something wrong. What if we continue our $\gamma$ so that it is more than just half of a great circle, or what if our curve $\gamma$ spirals many times around a point? Clearly, whenever $\gamma$ is not a piece of a great circle, the lunes will overlap. The correct definition for $|\gamma_\perp|$ is the area swept out by $\gamma(t)_\perp$, as $t$ varies, *counting multiplicities.* We now give the precise definitions.

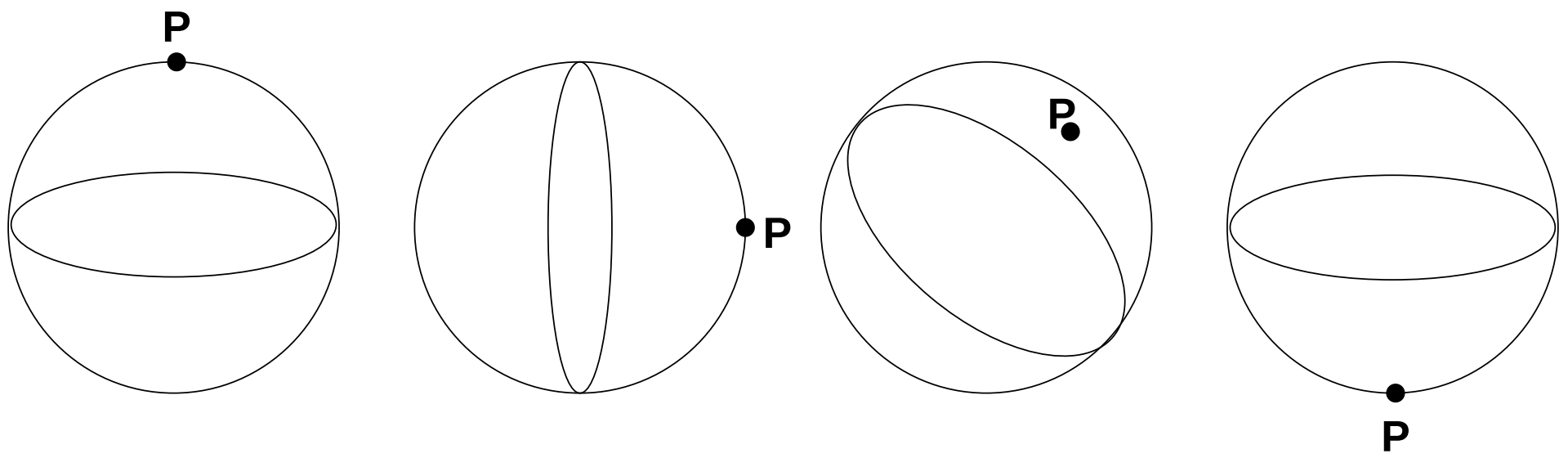

Figure 1. Points $P$ and associated equators $P_\perp$.

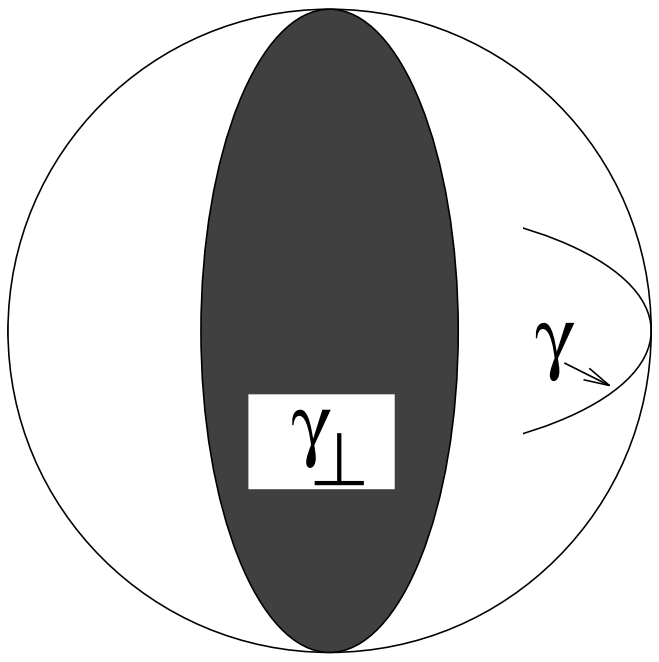

FIGURE 2. The lune $\cup\gamma_\perp$ when $\gamma$ is a great circular arc.

**Definition 2.3.** The *multiplicity* of a point $Q \in \cup\gamma_\perp$ is the number of equators in $\gamma_\perp$ that contain $Q$, i.e., the cardinality of $\{t \in \mathbb{R} | Q \in \gamma(t)_\perp\}$.

**Definition 2.4.** We define $|\gamma_\perp|$ to be the area of $\cup\gamma_\perp$ counting multiplicity. More precisely, we define $|\gamma_\perp|$ to be the integral of the multiplicity over $\cup\gamma_\perp$.

**Lemma 2.1.** *If $\gamma$ is a rectifiable curve, then*

$$\frac{|\gamma_\perp|}{\textit{area of } S^n} = \frac{|\gamma|}{\pi}. \tag{1}$$

As an example, consider a point $P$ on the surface of the Earth. If we assume that the point $P$ is receiving the direct ray of the sun—for our purposes, we consider the sun to be fixed in space relative to the Earth during the course of a day, with rays arriving in parallel—then $P_\perp$ is the great circle that divides day from night. This great circle is known to astronomers as the *terminator* (Figure 3). During the Earth's daily rotation, the point $P$ runs through all the points on a circle $\gamma$ of fixed latitude. Similarly, the Earth's rotation generates the collection of terminators $\gamma_\perp$.

The multiplicity in $\gamma_\perp$ is two on a region between two latitudes. This is a fancy mathematical way of saying that unless you are too close to the poles, you witness both a sunrise and a sunset every day! The summer solstice is a convenient example. $P$ is on the Tropic of Cancer and Equation (1) becomes

$$\frac{2 \times (\text{The surface area of the Earth between the Arctic/Antarctic Circles})}{\text{The surface area of the Earth}}$$
$$= \frac{\text{The length of the Tropic of Cancer}}{\pi \times (\text{The radius of the Earth})}$$

or equivalently

$$\frac{\text{The surface area of the Earth between the Arctic/Antarctic Circles}}{\text{The surface area of the Earth}}$$
$$= \frac{\text{The length of the Tropic of Cancer}}{\text{The length of the Equator}}.$$

2.2. **The expected number of real zeros of a random polynomial.** What does the geometric argument in the previous section and formula (1) in particular have to do with the number of real zeros of a random polynomial? Let

$$p(x) = a_0 + a_1 x + \cdots + a_n x^n$$

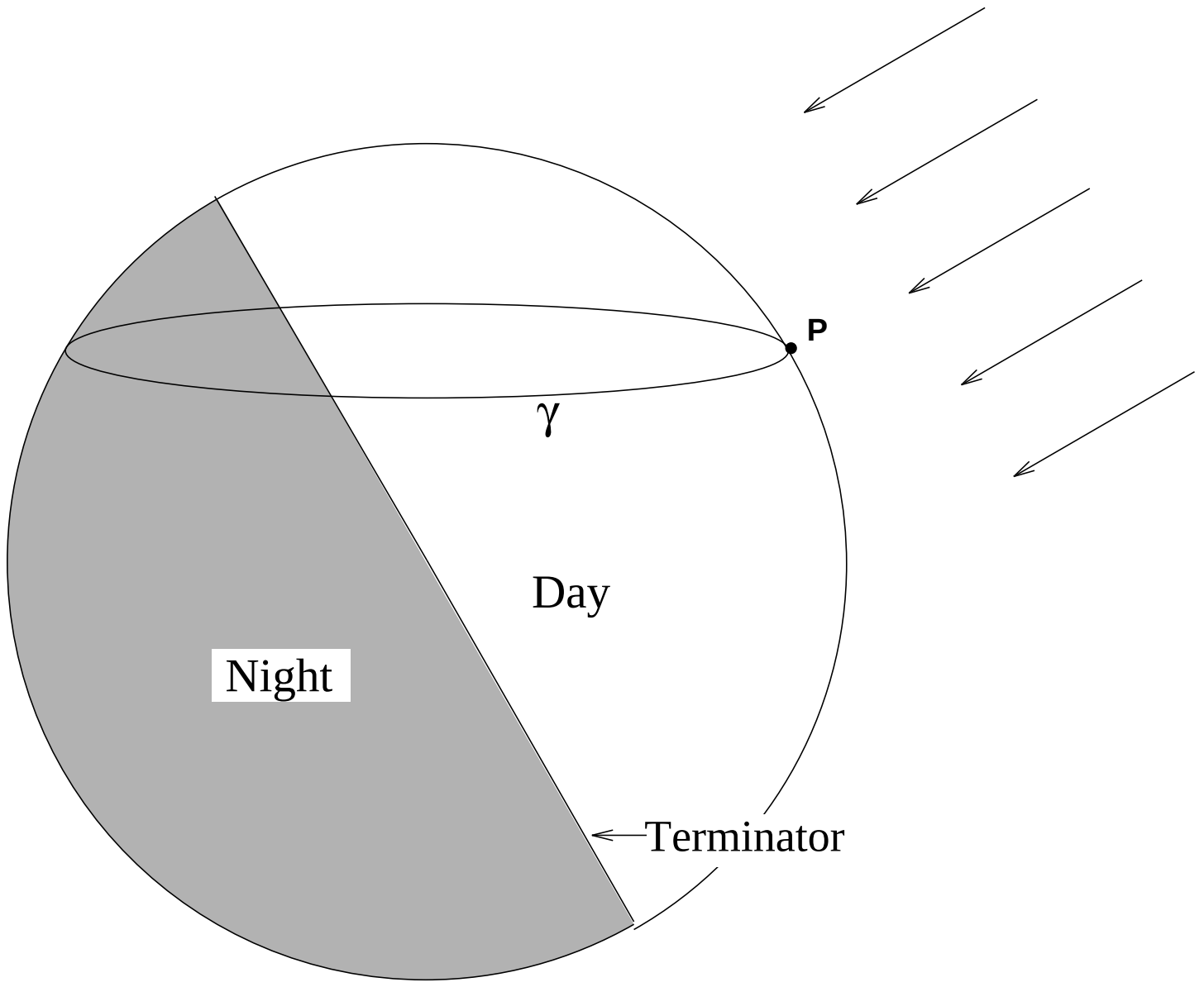

FIGURE 3. On the summer solstice, the direct ray of the sun reaches $P$ on the Tropic of Cancer $\gamma$.

be a non-zero polynomial. Define the two vectors

$$a = \begin{pmatrix} a_0 \\ a_1 \\ a_2 \\ \vdots \\ a_n \end{pmatrix} \text{ and } v(t) = \begin{pmatrix} 1 \\ t \\ t^2 \\ \vdots \\ t^n \end{pmatrix}.$$

The curve in $\mathbb{R}^{n+1}$ traced out by $v(t)$ as $t$ runs over the real line is called the *moment curve*.

The condition that $x = t$ is a zero of the polynomial $a_0 + a_1 x + \cdots + a_n x^n$ is precisely the condition that $a$ is perpendicular to $v(t)$. Another way of saying this is that $v(t)_\perp$ is the set of polynomials which have $t$ as a zero.

Define unit vectors

$$\mathbf{a} \equiv a/\|a\|, \qquad \gamma(t) \equiv v(t)/\|v(t)\|.$$

As before, $\gamma(t)_\perp$ corresponds to the polynomials which have $t$ as a zero.

When $n = 2$, the curve $\gamma$ is the intersection of an elliptical (squashed) cone and the unit sphere. In particular, $\gamma$ is not planar. If we include the point at infinity, $\gamma$ becomes a simple closed curve when $n$ is even. (In projective space, the curve is closed for all $n$.) The number of times that a point $\mathbf{a}$ on our sphere is covered by an equator is the multiplicity of $\mathbf{a}$ in $\gamma_\perp$. This is exactly the number of real zeros of the corresponding polynomial.

So far, we have not discussed *random* polynomials. If the $a_i$ are independent standard normals, then the vector $\mathbf{a}$ is uniformly distributed on the sphere $S^n$ since the joint density function in spherical coordinates is a function of the radius alone.

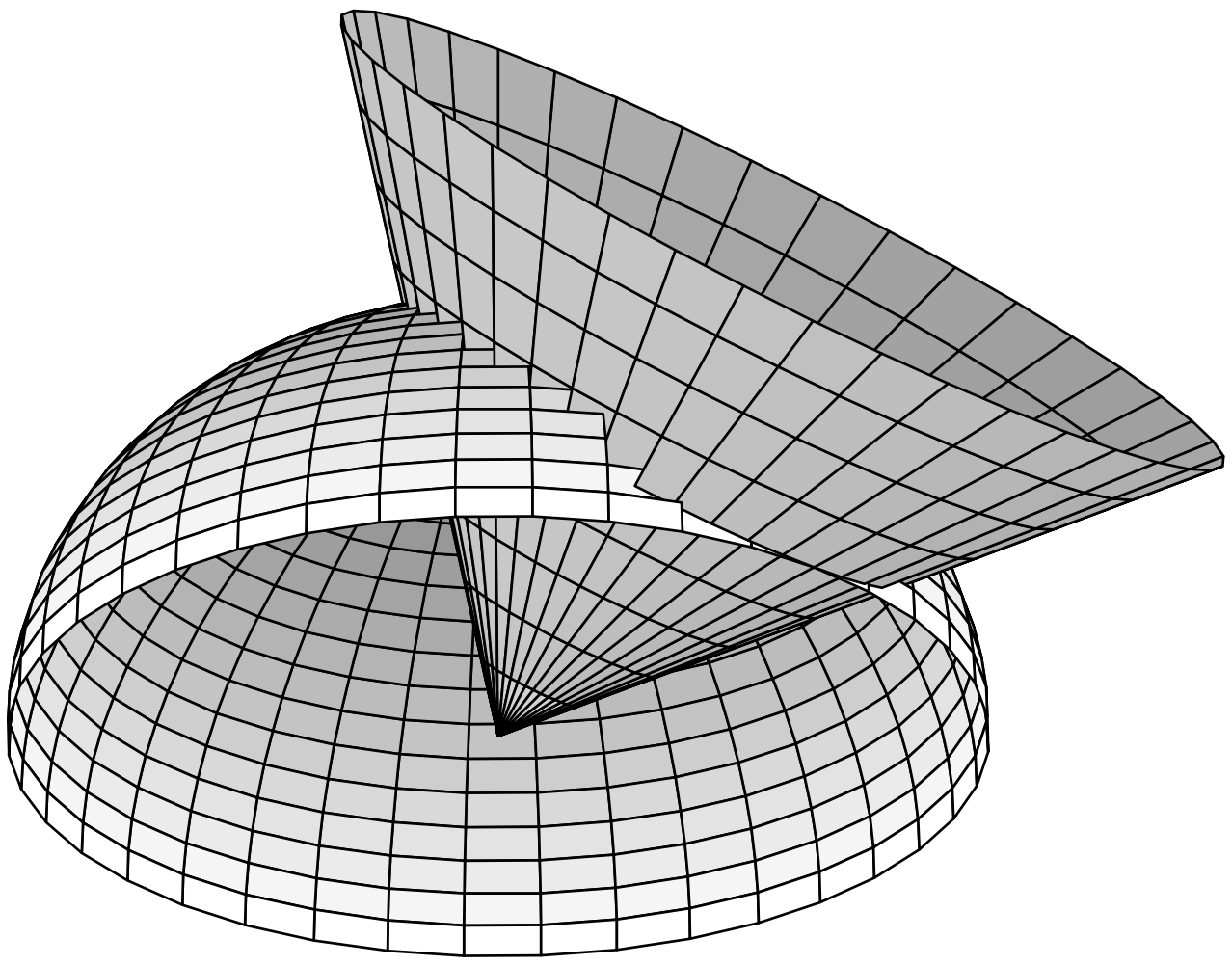

FIGURE 4. When $n = 2$, $\gamma$ is the intersection of the sphere and cone. The intersection is a curve that includes the North Pole and a point on the Equator.

What is $E_n \equiv$ the expected number of real zeros of a random polynomial? A random polynomial is identified with a uniformly distributed random point on the sphere, so $E_n$ is the area of the sphere with our convention of counting multiplicities.

Equation (1) (read backwards!) states that

$$E_n = \frac{1}{\pi}|\gamma|.$$

Our question about the expected number of real zeros of a random polynomial is reduced to finding the length of the curve $\gamma$. We compute this length in Section 2.3.

When $n = 2$, $\gamma$ is the intersection of the sphere and cone (Figure 4). The intersection is a curve that includes the North Pole and a point on the Equator.

2.3. **Calculating the length of $\gamma$.** We invoke calculus to obtain the integral formula for the length of $\gamma$ and hence the expected number of zeros of a random polynomial. The result was first obtained by Kac in 1943.

**Theorem 2.1** (Kac formula)**.** *The expected number of real zeros of a degree $n$ polynomial with independent standard normal coefficients is*

$$(2) \qquad \begin{aligned} E_n &= \frac{1}{\pi}\int_{-\infty}^{\infty}\sqrt{\frac{1}{(t^2-1)^2} - \frac{(n+1)^2 t^{2n}}{(t^{2n+2}-1)^2}}\,dt \\ &= \frac{4}{\pi}\int_0^1 \sqrt{\frac{1}{(1-t^2)^2} - \frac{(n+1)^2 t^{2n}}{(1-t^{2n+2})^2}}\,dt. \end{aligned}$$

*Proof.* The standard arclength formula is

$$|\gamma| = \int_{-\infty}^{\infty} \|\gamma'(t)\| \, dt.$$

We may proceed in two different ways.

*Method* I (Direct approach). To calculate the integrand, we first consider any differentiable $v(t) : \mathbb{R} \to \mathbb{R}^{n+1}$. It is not hard to show that

$$\gamma'(t) = \left( \frac{v(t)}{\sqrt{v(t) \cdot v(t)}} \right)' = \frac{[v(t) \cdot v(t)]v'(t) - [v(t) \cdot v'(t)]v(t)}{[v(t) \cdot v(t)]^{3/2}},$$

and therefore,

$$\begin{aligned} \|\gamma'(t)\|^2 &= \left( \frac{v(t)}{\sqrt{v(t) \cdot v(t)}} \right)' \cdot \left( \frac{v(t)}{\sqrt{v(t) \cdot v(t)}} \right)' \\ &= \frac{[v(t) \cdot v(t)][v'(t) \cdot v'(t)] - [v(t) \cdot v'(t)]^2}{[v(t) \cdot v(t)]^2}. \end{aligned}$$

If $v(t)$ is the moment curve, then we may calculate $\|\gamma'(t)\|$ with the help of the following observations and some messy algebra:

$$\begin{aligned} v(t) \cdot v(t) &= 1 + t^2 + t^4 + \cdots + t^{2n} = \frac{1 - t^{2n+2}}{1 - t^2}; \\ v'(t) \cdot v(t) &= t + 2t^3 + 3t^5 + \cdots + nt^{2n-1} \\ &= \frac{1}{2} \frac{d}{dt} \left( \frac{1 - t^{2n+2}}{1 - t^2} \right) = \frac{t \left(1 - t^{2\,n} - n\,t^{2\,n} + n\,t^{2\,n+2}\right)}{\left(t^2 - 1\right)^2}; \\ v'(t) \cdot v'(t) &= 1 + 4t^2 + 9t^4 + \cdots + n^2 t^{2n-2} \\ &= \frac{1}{4t} \frac{d}{dt} t \frac{d}{dt} \left( \frac{1 - t^{2n+2}}{1 - t^2} \right) = \frac{t^{2\,n+2} - t^2 - 1 + t^{2\,n} \left(n\,t^2 - n - 1\right)^2}{\left(t^2 - 1\right)^3}. \end{aligned}$$

Thus we arrive at the Kac formula:

$$\begin{aligned} E_n &= \frac{1}{\pi} \int_{-\infty}^{\infty} \frac{\sqrt{(t^{2n+2} - 1)^2 - (n+1)^2 t^{2n} (t^2 - 1)^2}}{(t^2 - 1)(t^{2n+2} - 1)} \, dt \\ &= \frac{1}{\pi} \int_{-\infty}^{\infty} \sqrt{\frac{1}{(t^2 - 1)^2} - \frac{(n+1)^2 t^{2n}}{(t^{2n+2} - 1)^2}} \, dt. \end{aligned}$$

*Method* II (Sneaky version). By introducing a logarithmic derivative, we can avoid the messy algebra in Method I. Let $v(t) : \mathbb{R} \to \mathbb{R}^{n+1}$ be any differentiable curve. Then it is easy to check that

$$\frac{\partial^2}{\partial x \partial y} \log[v(x) \cdot v(y)] \Big|_{y=x=t} = \|\gamma'(t)\|^2. \tag{3}$$

Thus we have an alternative expression for $\|\gamma'(t)\|^2$.

When $v(t)$ is the moment curve,

$$v(x) \cdot v(y) = 1 + xy + x^2 y^2 + \cdots + x^n y^n = \frac{1 - (xy)^{n+1}}{1 - xy},$$

the Kac formula is then

$$E_n = \frac{1}{\pi}\int_{-\infty}^{\infty} \sqrt{\frac{\partial^2}{\partial x \partial y} \log \frac{1-(xy)^{n+1}}{1-xy}\bigg|_{y=x=t}}\, dt.$$

This version of the Kac formula first appeared in [31]. In Section 4.4, we relate this sneaky approach to the so-called "Fubini-Study" metric. □

2.4. **The density of zeros.** Up until now, we have focused on the length of $\gamma = \{\gamma(t) | -\infty < t < \infty\}$ and concluded that it equals the expected number of zeros on the real line multiplied by $\pi$. What we really did, however, was compute the density of real zeros. Thus

$$\rho_n(t) \equiv \frac{1}{\pi}\sqrt{\frac{1}{(t^2-1)^2} - \frac{(n+1)^2 t^{2n}}{(t^{2n+2}-1)^2}}$$

is the expected number of real zeros per unit length at the point $t \in \mathbb{R}$. This is a true density: integrating $\rho_n(t)$ over any interval produces the expected number of real zeros on that interval. The *probability* density for a random real zero is $\rho_n(t)/E_n$. It is straightforward [26, 27] to see that as $n \to \infty$, the real zeros are concentrated near the point $t = \pm 1$.

The asymptotic behavior of both the density and expected number of real zeros is derived in the subsection below.

2.5. **The asymptotics of the Kac formula.** A short argument could have shown that $E_n \sim \frac{2}{\pi}\log n$ [26], but since several researchers, including Christensen, Sambandham, Stevens, and Wilkins have sharpened Kac's original estimate, we show here how successive terms of the asymptotic series may be derived, although we will derive only a few terms of the series explicitly. The constant $C_1$ and the next term $\frac{2}{n\pi}$ were unknown to previous researchers. See [2, pp. 90–91] for a summary of previous estimates of $C_1$.

**Theorem 2.2.** *As $n \to \infty$,*

$$E_n = \frac{2}{\pi}\log(n) + C_1 + \frac{2}{n\pi} + O(1/n^2)\ ,$$

*where*

$$C_1 = \frac{2}{\pi}\left(\log(2) \ + \int_0^\infty \left\{\sqrt{\frac{1}{x^2} - \frac{4e^{-2x}}{(1-e^{-2x})^2}} \ - \ \frac{1}{x+1}\right\} dx \right)$$
$$= 0.6257358072... \ .$$

*Proof.* We now study the asymptotic behavior of the density of zeros. To do this, we make the change of variables $t = 1 + x/n$, so

$$E_n = 4\int_0^\infty \hat{\rho}_n(x)\ dx\ ,$$

where

$$\hat{\rho}_n(x) = \frac{1}{n\pi}\sqrt{\frac{n^4}{x^2(2n+x)^2} - \frac{(n+1)^2(1+x/n)^{2n}}{[(1+x/n)^{2n+2}-1]^2}}$$

is the (transformed) density of zeros. Using

$$\left(1+\frac{x}{n}\right)^n = e^x\left(1-\frac{x^2}{2n}\right) + O(1/n^2)\ ,$$

we see that for any fixed $x$, as $n \to \infty$, the density of zeros is given by

$$\hat{\rho}_n(x) = \hat{\rho}_\infty(x) + \left[\frac{x(2-x)}{2n}\hat{\rho}_\infty(x)\right]' \; + \; O(1/n^2) \; , \tag{4}$$

where

$$\hat{\rho}_\infty(x) \equiv \frac{1}{2\pi}\left[\frac{1}{x^2} - \frac{4e^{-2x}}{(1-e^{-2x})^2}\right]^{1/2} \; .$$

This asymptotic series cannot be integrated term by term. We solve this problem by noting that

$$\frac{\chi[x>1]}{2\pi x} - \frac{1}{2\pi(2n+x)} = \frac{\chi[x>1]}{2\pi x} - \frac{1}{4n\pi} \; + \; O(1/n^2) \; , \tag{5}$$

where we have introduced the factor

$$\chi[x>1] \equiv \begin{cases} 1 & \text{if } x > 1, \\ 0 & \text{if } x \le 1 \end{cases}$$

to avoid the pole at $x = 0$. Subtracting (5) from (4), we obtain

$$\begin{aligned} &\hat{\rho}_n(x) - \left\{\frac{\chi[x>1]}{2\pi x} - \frac{1}{2\pi(2n+x)}\right\} \\ &\quad = \left\{\hat{\rho}_\infty(x) - \frac{\chi[x>1]}{2\pi x}\right\} + \left\{\left[\frac{x(2-x)}{2n}\hat{\rho}_\infty(x)\right]' + \frac{1}{4\pi n}\right\} \; + \; O(1/n^2) \; . \end{aligned}$$

We then integrate term by term from 0 to $\infty$ to get

$$\begin{aligned} &\int_0^\infty \hat{\rho}_n(x)\; dx \; - \; \frac{1}{2\pi}\log(2n) \\ &\quad = \; \int_0^\infty \left\{\hat{\rho}_\infty(x) - \frac{\chi[x>1]}{2\pi x}\right\} dx + \frac{1}{2n\pi} \; + \; O(1/n^2) \; . \end{aligned}$$

The theorem immediately follows from this formula and one final trick: we replace $\chi[x>1]/x$ with $1/(x+1)$ in the definition of $C_1$ so we can express it as a single integral of an *elementary* function. □

## 3. Random functions with central normal coefficients

Reviewing the discussion in Section 2, we see that we could omit some members of our basis set $\{1, x, x^2, \ldots, x^n\}$ and ask how many zeros are expected to be real of an $n$th degree polynomial with, say, its cubic term deleted. The proof would hardly change. Or we can change the function space entirely and ask how many zeros of the random function

$$a_0 + a_1 \sin(x) + a_2 e^{|x|}$$

are expected to be real—the answer is 0.63662. The only assumption is that the coefficients are independent standard normals. If $f_0, f_1, \ldots, f_n$ is any collection of rectifiable functions, we may define the analogue of the moment curve

$$v(t) = \begin{pmatrix} f_0(t) \\ f_1(t) \\ \vdots \\ f_n(t) \end{pmatrix} . \tag{6}$$

The function $\frac{1}{\pi}\|\gamma'(t)\|$ is the density of a real zero; its integral over $\mathbb{R}$ is the expected number of real zeros.

We may relax the assumption that the coefficient vector $a = (a_0, \ldots, a_n)^T$ contains independent standard normals by allowing for any multivariate distribution with zero mean. If the $a_i$ are normally distributed, $E(a) = 0$ and $E(aa^T) = C$, then $a$ is a (central) multivariate normal distribution with covariance matrix $C$. It is easy to see that $a$ has this distribution if and only if $C^{-1/2}a$ is a vector of standard normals. Since

$$a \cdot v(t) = C^{-1/2}a \cdot C^{1/2}v(t),$$

the density of real zeros with coefficients from an arbitrary central multivariate normal distribution is

$$\frac{1}{\pi}\|\mathbf{w}'(t)\|, \text{ where } w(t) = C^{1/2}v(t), \quad \text{and} \quad \mathbf{w}(t) = w(t)/\|w(t)\|. \tag{7}$$

The expected number of real zeros is the integral of $\frac{1}{\pi}\|\mathbf{w}'(t)\|$.

We now state our general result.

**Theorem 3.1.** *Let $v(t) = (f_0(t), \ldots, f_n(t))^T$ be any collection of differentiable functions and $a_0, \ldots, a_n$ be the elements of a multivariate normal distribution with mean zero and covariance matrix $C$. The expected number of real zeros on an interval (or measurable set) $I$ of the equation*

$$a_0 f_0(t) + a_1 f_1(t) + \cdots + a_n f_n(t) = 0$$

*is*

$$\int_I \frac{1}{\pi}\|\mathbf{w}'(t)\| dt,$$

*where $\mathbf{w}$ is defined by Equations (7). In logarithmic derivative notation this is*

$$\frac{1}{\pi}\int_I \left(\frac{\partial^2}{\partial x \partial y} \left(\log v(x)^T C v(y)\right)\big|_{y=x=t}\right)^{1/2} dt.$$

Geometrically, changing the covariance is the same as changing the inner product on the space of functions.

We now enumerate several examples of Theorem 3.1. We consider examples for which $v(x)^T C v(y)$ is a nice enough function of $x$ and $y$ that the density of zeros can be easily described. For a survey of the literature, see [2], which also includes the results of numerical experiments. In our discussion of random series, proofs of convergence are omitted. Interested readers may refer to [45]. We also suggest the classic book of J.-P. Kahane [28], where other problems about random series of functions are considered.

## 3.1. **Random polynomials.**

3.1.1. *The Kac formula.* If the coefficients of random polynomials are independent standard normal random variables, we saw in the previous section that from

$$v(x)^T C v(y) = \frac{1-(xy)^{n+1}}{1-xy}, \tag{8}$$

we can derive the Kac formula.

3.1.2. *A random polynomial with a simple answer.* Consider random polynomials

$$a_0 + a_1 x + \cdots + a_n x^n,$$

where the $a_i$ are independent normals with variances $\binom{n}{i}$. Such random polynomials have been studied because of their mathematical properties [31, 46] and because of their relationship to quantum physics [4].

By the binomial theorem,

$$v(x)^T C v(y) = \sum_{k=0}^{n} \binom{n}{k} x^k y^k = (1+xy)^n.$$

We see that the density of zeros is given by

$$\rho(t) = \frac{\sqrt{n}}{\pi(1+t^2)}.$$

This is a Cauchy distribution, that is, $\arctan(t)$ is uniformly distributed on $[-\pi/2, \pi/2]$. Integrating the density shows that the expected number of real zeros is $\sqrt{n}$. As we shall see in Section 4.1, this simple expected value and density is reflected in the geometry of $\gamma$.

As an application, assume that $p(t)$ and $q(t)$ are independent random polynomials of degree $n$ with coefficients distributed as in this example. By considering the equation $p(t) - tq(t) = 0$, it is possible to show that the expected number of *fixed points* of the rational mapping

$$p(t)/q(t) : \mathbb{R} \cup \{\infty\} \to \mathbb{R} \cup \{\infty\}$$

is exactly $\sqrt{n+1}$.

3.1.3. *Application: Spijker's lemma on the Riemann sphere.* Any curve in $\mathbb{R}^n$ can be interpreted as $v(t)$ for some space of random functions. Let

$$r(t) = \frac{a(t) + ib(t)}{c(t) + id(t)}$$

be any rational function, where $a, b, c,$ and $d$ are real polynomials of a real variable $t$. Let $\gamma$ be the stereographic projection of $r(t)$ onto the Riemann sphere. It is not difficult to show that $\gamma$ is the projection of the curve

$$(f_0(t), f_1(t), f_2(t))$$

onto the unit (Riemann) sphere, where $f_0 = 2(ac+bd), f_1 = 2(bc-ad)$, and $f_2 = a^2 + b^2 - c^2 - d^2$. The geometry is illustrated in Figure 5.

Therefore the length of $\gamma$ is $\pi$ times the expected number of real zeros of the random function

$$a_0 f_0 + a_1 f_1 + a_2 f_2,$$

where the $a_i$ are independent standard normals. For example, if $a, b, c,$ and $d$ are polynomials of degrees no more than $n$, then any such function has degree at most $2n$, so the length of $\gamma$ can be no more than $2n\pi$. By taking a Möbius transformation, we arrive at Spijker's lemma:

*The image, on the Riemann sphere, of any circle under a complex rational mapping, with numerator and denominator having degrees no more than $n$, has length no longer than $2n\pi$.*

This example was obtained from Wegert and Trefethen [51].

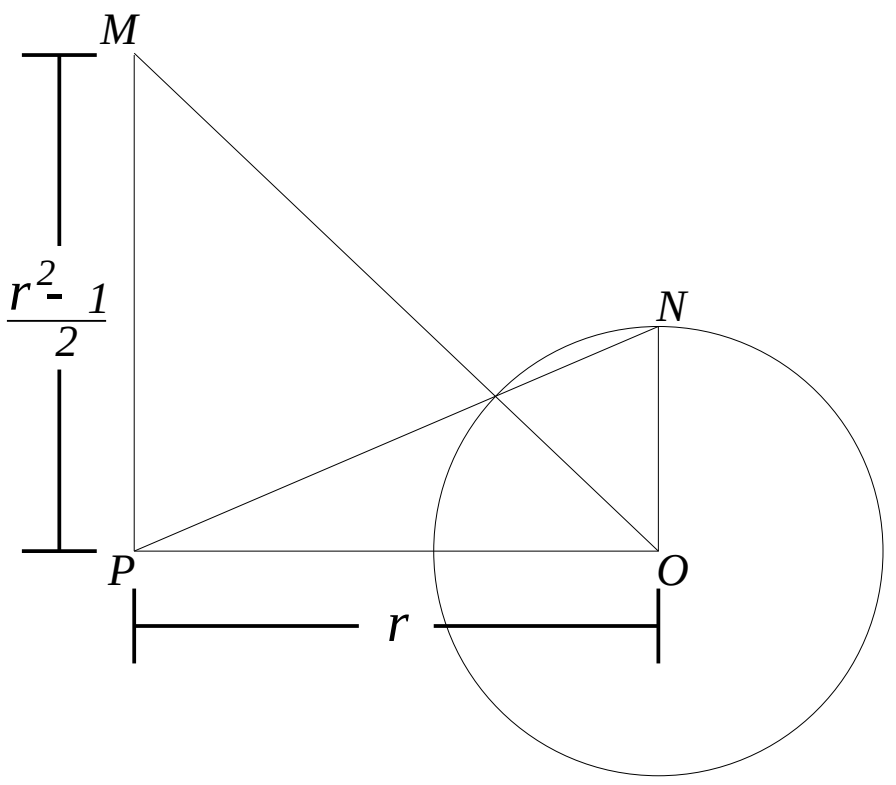

FIGURE 5

3.1.4. *Random sums of orthogonal polynomials.* Consider the vector space of polynomials of the form $\sum_{k=0}^{n} a_k P_k(x)$ where $a_k$ are independent standard normal random variables and where $\{P_k(x)\}$ is a set of normalized orthogonal polynomials with any non-negative weight function on any interval. The Darboux-Christoffel formula [21, 8.902] states that

$$\sum_{k=0}^{n} P_k(x)P_k(y) = \left(\frac{q_n}{q_{n+1}}\right) \frac{P_n(y)P_{n+1}(x) - P_n(x)P_{n+1}(y)}{x-y},$$

where $q_n$ (resp. $q_{n+1}$) is the leading coefficient of $P_n$ (resp. $P_{n+1}$). With this formula and a bit of work, we see that

$$\rho(t) = \frac{\sqrt{3}}{6\pi}\sqrt{2G'(t) - G^2(t)},$$

where

$$G(t) = \frac{d}{dt} \log \frac{d}{dt}\left(\frac{P_{n+1}(t)}{P_n(t)}\right).$$

This is equivalent to formula (5.21) in [2]. Interesting asymptotic results have been derived by Das and Bhatt. The easiest example to consider is that of random sums of Chebyshev polynomials, for which the density of zeros is an elementary function of $n$ and $t$.

3.2. **Random infinite series.**

3.2.1. *Power series with uncorrelated coefficients.* Consider a random power series

$$f(x) = a_0 + a_1 x + a_2 x^2 + \cdots,$$

where $a_k$ are independent standard normal random variables. This has radius of convergence one with probability one. Thus we will assume that $-1 < x < 1$. In this case,

$$v(x)^T C v(y) = \frac{1}{1-xy}.$$

The logarithmic derivative reveals a density of zeros of the form

$$\rho(t) = \frac{1}{\pi(1-t^2)}.$$

We see that the expected number of zeros on any subinterval $[a,b]$ of $(-1,1)$ is

$$\frac{1}{2\pi}\log\frac{(1-a)(1+b)}{(1+a)(1-b)}.$$

This result may also be derived from the original Kac formula by letting $n\to\infty$.

3.2.2. *Power series with correlated coefficients.* What effect does correlation have on the density of zeros? We will consider a simple generalization of the previous example. Consider the random power series

$$f(x)=a_0+a_1x+a_2x^2+\cdots,$$

where $a_k$ are standard normal random variables, but assume that the correlation between $a_k$ and $a_{k+1}$ equals some constant $r$ for all $k$. Thus the covariance matrix is tridiagonal with one on the diagonal and $r$ on the superdiagonal and subdiagonal. In order to assure that this matrix be positive definite, we will assume that $|r|\le\frac{1}{2}$. By the Gershgorin Theorem the spectral radius of the covariance matrix is less than or equal to $1+2r$, and therefore the radius of convergence of the random sequence is independent of $r$. Thus we will, as in the previous example, assume that $-1<x<1$. We see that

$$v(x)^TCv(y)=\frac{1+r(x+y)}{1-xy},$$

so

$$\rho(t)=\frac{1}{\pi}\sqrt{\frac{1}{(1-t^2)^2}-\frac{r^2}{(1+2rt)^2}}.$$

Notice that the correlation between coefficients has decreased the density of zeros throughout the interval.

3.2.3. *Random entire functions.* Consider a random power series

$$f(x)=a_0+a_1x+a_2x^2+\cdots,$$

where $a_k$ are independent central normal random variables with variances $1/k!$, i.e., the covariance matrix is diagonal with the numbers $1/k!$ down the diagonal. This series has infinite radius of convergence with probability one. Now clearly

$$v(x)^TCv(y)=e^{xy},$$

so $\rho(t)=1/\pi$. In other words, the real zeros are uniformly distributed on the real line, with a density of $1/\pi$ zeros per unit length.

3.2.4. *Random trigonometric sums and Fourier series.* Consider the trigonometric sum

$$\sum_{k=0}^{\infty}a_k\cos\nu_k\theta+b_k\sin\nu_k\theta,$$

where $a_k$ and $b_k$ are independent normal random variables with means zero and variances $\sigma_k^2$. Notice that

$$v(x)^TCv(y)=\sum_{k=0}^{\infty}\sigma_k^2(\sin\nu_kx\sin\nu_ky+\cos\nu_kx\cos\nu_ky)=\sum_{k=0}^{\infty}\sigma_k^2\cos\nu_k(x-y),$$

and we see that the density of roots is constant. Thus the real zeros of the random trigonometric sum are uniformly distributed on the real line, and the expected number of zeros on the interval $[a, b]$ is

$$\frac{b-a}{\pi}\sqrt{\frac{\sum \nu_k^2 \sigma_k^2}{\sum \sigma_k^2}}.$$

Note that the slower the rate of convergence of the series, the higher the root density. For example, if $\sigma_k = k^{-3/2}$ and $v_k = k$, then the series converges uniformly with probability one, but the root density is infinite.

The similarity between this formula and the Pythagorean theorem is more than superficial, as we will see when we discuss the geodesics of flat tori in Section 4.2. Several authors, including Christensen, Das, Dunnage, Jamrom, Maruthachalam, Qualls, and Sambandham [2] have derived results about the expected number of zeros of these and other trigonometric sums.

3.2.5. *Random Dirichlet series.* Consider a random Dirichlet series

$$f(x) = a_1 + \frac{a_2}{2^x} + \frac{a_3}{3^x} + \cdots,$$

where $a_k$ are independent standard normal random variables. This converges with probability one if $x > 1/2$. We see that

$$v(x)^T C v(y) = \sum_{k=1}^{\infty} \frac{1}{k^{x+y}} = \zeta(x+y)$$

and that the expected number of zeros on any interval $[a, b], a > 1/2$, is

$$\frac{1}{2\pi}\int_a^b \sqrt{[\log \zeta(2t)]''}\, dt.$$

## 4. Theoretical considerations

4.1. **A curve with more symmetries than meet the eye.** We return to the example in Section 3.1.2 and explore why the distribution of the real zeros was so simple. Take the curve $v(t)$ and make the change of variables $t = \tan\theta$ and scale to obtain

$$\gamma(\theta) = (\cos^n\theta)(v(\tan\theta)).$$

Doing so shows that

$$\gamma(\theta) = \begin{pmatrix} \binom{n}{0}^{1/2} \cos^n\theta \\ \binom{n}{1}^{1/2} \cos^{n-1}\theta\sin\theta \\ \binom{n}{2}^{1/2} \cos^{n-2}\theta\sin^2\theta \\ \vdots \\ \binom{n}{n}^{1/2} \sin^n\theta \end{pmatrix},$$

i.e., $\gamma_k(\theta) = \binom{n}{k}^{1/2} \cos^{n-k}\theta \sin^k\theta$, where the dimension index $k$ runs from 0 to $n$.

We have chosen to denote this curve the *super-circle*. The binomial expansion of $(\cos^2\theta + \sin^2\theta)^n = 1$ tells us that our super-circle lives on the unit sphere. Indeed when $n = 2$, the super-circle is merely a small-circle on the unit sphere in $\mathbb{R}^3$. When $n = 1$, the super-circle is the unit circle in the plane.

What is not immediately obvious is that every point on this curve "looks" exactly the same. We will display an orthogonal matrix $Q(\phi)$ that rotates $\mathbb{R}^{n+1}$ in such a manner that each and every point on the super-circle $\gamma(\theta)$ is sent to $\gamma(\theta+\phi)$. To do this, we show that $\gamma$ is a solution to a "nice" ordinary differential equation.

By a simple differentiation of the $k$th component of $\gamma(\theta)$, we see that

$$\frac{d}{d\theta}\gamma_k(\theta) = \alpha_k \gamma_{k-1}(\theta) - \alpha_{k+1}\gamma_{k+1}(\theta), \quad k = 0, \dots, n,$$

where $\alpha_k \equiv \sqrt{k(n+1-k)}$. In matrix–vector notation this means that

(9)

$$\frac{d}{d\theta}\gamma(\theta) = A\gamma(\theta), \text{ where } A = \begin{pmatrix} 0 & -\alpha_1 & & & & \\ \alpha_1 & 0 & -\alpha_2 & & & \\ & \alpha_2 & 0 & -\alpha_3 & & \\ & & \ddots & \ddots & \ddots & \\ & & & \alpha_{n-1} & 0 & -\alpha_n \\ & & & & \alpha_n & 0 \end{pmatrix};$$

i.e., $A$ has the $\alpha_i$ on the subdiagonal, the $-\alpha_i$ on the superdiagonal, and 0 everywhere else, including the main diagonal.

The solution to the ordinary differential equation (9) is

$$\gamma(\theta) = e^{A\theta}\gamma(0). \tag{10}$$

The matrix $Q(\phi) \equiv e^{A\phi}$ is orthogonal because $A$ is anti-symmetric, and indeed $Q(\phi)$ is the orthogonal matrix that we promised would send $\gamma(\theta)$ to $\gamma(\theta+\phi)$. We suspect that (10) with the specification that $\gamma(0) = (1, 0, \dots, 0)^T$ is the most convenient description of the super-circle. Differentiating (10) any number of times shows explicitly that

$$\frac{d^j\gamma}{d\theta^j}(\theta) = e^{A\theta}\frac{d^j\gamma}{d\theta^j}(0).$$

In particular, the speed is invariant. A quick check shows that it is $\sqrt{n}$. If we let $\theta$ run from $-\pi/2$ to $\pi/2$, we trace out a curve of length $\pi\sqrt{n}$.

The ideas here may also be expressed in the language of invariant measures for polynomials [31]. This gives a deeper understanding of the symmetries that we will only sketch here. Rather than representing a polynomial as

$$p(t) = a_0 + a_1 t + a_2 t^2 + \cdots + a_n t^n,$$

we homogenize the polynomial and consider

$$\hat{p}(t_1, t_2) = a_0 t_2^n + a_1 t_1 t_2^{n-1} + \cdots + a_{n-1} t_1^{n-1} + a_n t_1^n.$$

For any angle $\alpha$, a new "rotated" polynomial may be defined by

$$\hat{p}_\alpha(t_1, t_2) = \hat{p}(t_1\cos\alpha + t_2\sin\alpha, -t_1\sin\alpha + t_2\cos\alpha).$$

It is not difficult to show directly that if the $a_i$ are independent and normally distributed with variance $\binom{n}{i}$, then so are the coefficients of the rotated polynomial. The symmetry of the curve and the symmetry of the polynomial distribution are equivalent. An immediate consequence of the rotational invariance is that the distribution of the real zeros must be Cauchy.

4.2. **Geodesics on flat tori.** We now a take a closer look at the random trigonometric sums in Section 3.2.4. Fix a finite interval $[a, b]$. For simplicity assume that

$$\sum_{k=0}^{n} \sigma_k^2 = 1 \ .$$

The curve $\gamma(\theta)$ is given by

$$(\sigma_0 \cos \nu_0 \theta, \sigma_0 \sin \nu_0 \theta, \ldots, \sigma_n \cos \nu_n \theta, \sigma_n \sin \nu_n \theta) \ .$$

This curve is a geodesic on the flat $(n+1)$-dimensional torus

$$(\sigma_0 \cos \theta_0, \sigma_0 \sin \theta_0, \ldots, \sigma_n \cos \theta_n, \sigma_n \sin \theta_n) \ .$$

Therefore if we lift to the universal covering space of the torus, $\gamma$ becomes a straight line in $\mathbb{R}^{n+1}$. By the Pythagorean theorem, the length of $\gamma$ is

$$(b-a)\sqrt{\sum_{k=0}^{n} \nu_k^2 \sigma_k^2} \ ,$$

which equals $\pi$ times the expected number of zeros on the interval $[a, b]$.

Now replace $[a, b]$ with $(-\infty, \infty)$. If $\nu_i/\nu_j$ is rational for all $i$ and $j$, then $\gamma$ is closed; otherwise it is dense in some subtorus.

Now consider the $\gamma(\theta)$ discussed in Section 4.1. Observe that $\gamma(x)^T \gamma(y) = \cos^n(x-y)$. Thus *if we choose* the $\nu_k$ and the $\sigma_k$ correctly, the polynomial example in Section 3.1.2 becomes a special case of a random trigonometric sum. Thus the super-circle discussed in Section 4.1 is a geodesic on a flat torus.

4.3. **The Kac matrix.** Mark Kac was the first mathematician to obtain an exact formula for the expected number of real zeros of a random polynomial. Ironically, he also has his name attached to a certain matrix that is important to understanding random polynomials, yet we have no evidence that he ever made the connection.

The $(n+1) \times (n+1)$ Kac matrix is defined as the tridiagonal matrix

$$S_n = \begin{pmatrix} 0 & n & & & & \\ 1 & 0 & n-1 & & & \\ & 2 & 0 & n-2 & & \\ & & \ddots & \ddots & \ddots & \\ & & & n-1 & 0 & 1 \\ & & & & n & 0 \end{pmatrix} .$$

The history of this matrix is documented in [49], where there are several proofs that $S_n$ has eigenvalues $-n, -n+2, -n+4, \ldots, n-2, n$. One of the proofs is denoted as "mild trickery by Kac". We will derive the eigenvalues by employing a different trick.

**Theorem 4.1.** *The eigenvalues of $S_n$ are the integers $2k - n$ for $k = 0, 1, \ldots, n$.*

*Proof.* Define

$$f_k(x) \equiv \sinh^k(x) \cosh^{n-k}(x), \quad k = 0, \ldots, n,$$
$$g_k(x) \equiv (\sinh(x) + \cosh(x))^k (\sinh(x) - \cosh(x))^{n-k}, \quad k = 0, \ldots, n.$$

If $V$ is the vector space of functions with basis $\{f_k(x)\}$, then the $g_k(x)$ are clearly in this vector space. Also, $\frac{d}{dx} f_k(x) = k f_{k-1}(x) + (n-k) f_{k+1}(x)$, so that the Kac matrix is the representation of the operator $d/dx$ in $V$. We actually wrote $g_k(x)$ in a more complicated way than we needed to so that we could emphasize that $g_k(x) \in V$. Actually, $g_k(x) = \exp((2k-n)x)$ is an eigenfunction of $d/dx$ with eigenvalue $2k-n$ for $k = 0, \ldots, n$. The eigenvector is obtained by expanding the above expression for $g_k(x)$. □

An alternative tricky proof using graph theory is to consider the $2^n \times 2^n$ incidence matrix of an $n$-dimensional hypercube. This matrix is the tensor (or Kronecker) product of $\left(\begin{smallmatrix} 0 & 1 \\ 1 & 0 \end{smallmatrix}\right)$ $n$ times, so the eigenvalues of this matrix are sums of the form $\sum_{i=1}^{n} \pm 1$, i.e., this matrix has $2^n$ eigenvalues all of which have the form $2k-n$ for $k = 0, \ldots, n$. This matrix is closely related to the discrete Laplacian of the hypercube graph and the $n$-fold discrete Fourier transform on a grid with edge length 2. So far we have the right set of eigenvalues but the wrong matrix. However, if we collapse the matrix by identifying those nodes with $k = 0, 1, \ldots, n$ ones in their binary representation, we obtain the $(n+1) \times (n+1)$ Kac matrix transposed. (Any node with $k$ ones has $k$ neighbors with $k-1$ ones and $n-k$ neighbors with $k+1$ ones.) It is an interesting exercise to check that by summing eigenvectors over all possible symmetries of the hypercube, the projected operator inherits the eigenvalues $2k-n$ ($k = 0, \ldots, n$), each with multiplicity 1. □

We learned of this second tricky proof from Persi Diaconis, who explained it to us in terms of random walks on the hypercube and the Ehrenfest urn model of diffusion [8, 9]. The Kac matrix is also known as the "Clement matrix" in Higham's *Test matrix toolbox for Matlab* [25] because of Clement's [7] proposed use of this matrix as a test matrix. Numerically, it can be quite difficult to obtain all of these integer eigenvalues.

The symmetrized Kac matrix looks exactly like the matrix $A$ in (9) without any minus signs. Indeed $iS_n$ is similar to the matrix in (9).

4.4. **The Fubini-Study metric.** We now reveal the secret that inspired the "sneaky" approach to the calculation of the length of the curve $\gamma(t) = v(t)/\|v(t)\|$ that appears in Section 2.3. (See (3).) The secret that we will describe is the Fubini-Study metric.

An interesting struggle occurs in mathematics when quotient spaces are defined. Psychologically, it is often easier to think of an individual representative of an equivalence class rather than the class itself. As mathematicians, we train ourselves to overcome this; but practically speaking, when it is time to compute, we still must choose a representative. As an example, consider vectors $v \in \mathbb{R}^n$ and its projection $v/\|v\|$ onto the sphere. (If we do not distinguish $\pm v/\|v\|$, we are then in projective space.) The normalization obtained from the division by $\|v\|$ is a distraction that we would like to avoid.

Perhaps a more compelling example may be taken from the set of $n \times p$ matrices $M$ with $n > p$. The Grassman manifold is obtained by forming the equivalence class of all rank $p$ matrices $M$ whose columns span the same subspace of $\mathbb{R}^n$. To compute a canonical form for $M$ may be an unnecessary bother that we would like to avoid. When $p = 1$, the Grassman manifold reduces to the projective space example in the previous paragraph.

The Fubini-Study metric on projective space allows us to keep the $v$ for our coordinates in the first example. The more general version for the Grassman manifold allows us to keep the $M$. A historical discussion of Fubini's original ideas may be found in [35]. We have seen only the complex version in the standard texts [22, 29], but for simplicity, we discuss the real case here.

We see that $\gamma(t)$ is independent of $\|v(t)\|$, i.e., it is invariant under scaling. The logarithmic derivative is specifically tailored to be invariant under scaling by any $\lambda(t)$:

$$\frac{\partial^2}{\partial x \partial y} \log[\lambda(x)v(x) \cdot \lambda(y)v(y)] = \frac{\partial^2}{\partial x \partial y} \{\log[v(x) \cdot v(y)] + \log(\lambda(x)) + \log(\lambda(y))\}$$
$$= \frac{\partial^2}{\partial x \partial y} \log[v(x) \cdot v(y)].$$

The logarithmic derivative may appear complicated, but it is a fair price to pay to eliminate $\|v(t)\|$. The length of the projected version of $v(t)$ traced out by $t \in [a, b]$ is

$$\int_a^b \sqrt{\frac{\partial^2}{\partial x \partial y} \log[v(x) \cdot v(y)]\bigg|_{y=x=t}}\, dt.$$

The integrand is the square root of the determinant of the metric tensor. This is the "pull-back" of a metric tensor on projective space.

The Grassman version is almost the same; it takes the form

$$\int_a^b \sqrt{\frac{\partial^2}{\partial x \partial y} \log\det[M(x)^{\mathrm{T}} M(y)],\bigg|_{y=x=t}}\, dt,$$

where $M(t)$ denotes a curve in matrix space.

4.5. **Integral geometry.** Integral geometry (sometimes known as Geometric Probability) relates the measures of random manifolds and their intersections. References (such as [43], [47], [44, p. 253], and [5, p. 73]) vary in terms of setting and degree of generality.

For our purposes we will consider two submanifolds $M$ and $N$ of the sphere $S^{m+n}$, where $M$ has dimension $m$ and $N$ has dimension $n$. If $Q$ is a random orthogonal matrix (i.e., a random rotation), then

$$E(\#(M \cap QN)) = \frac{2}{|S^m||S^n|} |M||N|. \tag{11}$$

In words, the formula states that the expected number of intersections of $M$ with a randomly rotated $N$ is twice the product of the volumes of $M$ and $N$ divided by the product of the volumes of spheres.

For us "number of intersections" has the interpretation of "number of zeros" so that we may relate the average number of zeros with the lengths of curves (or more generally volumes of surfaces). We will apply this formula directly when we consider random systems of equations in Section 7.

If the manifold $N$ is itself random and independent of $Q$, then the formula above is correct with the understanding that $|N|$ refers to the average volume of $N$. This formulation is needed for Lemma 6.1.

The factor of 2 often disappears in practical computations. Mathematically, all of the action is on the half-sized projective space rather than on the sphere.

4.6. **The evaluation mapping.** The defining property of a *function* space is that its elements can be evaluated. To be precise, if $F$ is a vector space of real-valued functions defined on some set $S$, we have an evaluation mapping, $ev : S \to F^*$, defined by $ev(s)(f) = f(s)$, that tells us everything about the function space. Conversely, if we are given *any* vector space $F$ and *any* function from $S$ to $F^*$, we may take this function to be the evaluation mapping and thus convert $F$ into a function space.

Pick an element $f$ of $F$ (at random). The annihilator of $f$ is the set $f_\perp = \{\theta \in F^* | \theta(f) = 0\}$. Checking definitions, we see that the intersections of $f_\perp$ with the image of $ev$ correspond to zeros of $f$. Thus the average number of intersections is the average number of zeros.

Now let us choose an inner product for $F$, or equivalently, let us choose a central normal measure for $F$. If $ev(S)$ is a rectifiable curve, we may apply integral geometry and conclude the following:

**Theorem 4.2.** *The expected number of zeros is the length of the projection of the image of the evaluation mapping onto the unit sphere in the dual space divided by* $\pi$.

Thus the expected number of zeros is proportional to the "size" of the image of the evaluation mapping.

The inner product also gives rise to an isomorphism $\iota : F \to F^*$, defined by $\iota(f)(g) = f \cdot g$. It is just a matter of checking definitions to see that

$$v(t) = \iota^{-1} ev(t)$$

is *the dual of the evaluation mapping.* Thus $v(t)$ is the natural object that describes both the function space $F$ and the choice of inner product.

## 5. Extensions to other distributions

This paper began by considering random polynomials with standard normal coefficients, and then we realized quickly that any multivariate normal distribution with mean zero (the so-called "central distributions") hardly presented any further difficulty. We now generalize to arbitrary distributions, with a particular focus on the non-central multivariate normal distributions. The basic theme is the same: the density of zeros is equal to the rate at which the equators of a curve sweeps out area. Previous investigations are surveyed in [2]. In the closely related work of Rice [41] and [42, p. 52], expressions are obtained for the distributions of zeros. Unfortunately, these expressions appeared unwieldy for computing even the distribution for the quadratic [41, p. 414]. There is also the interesting recent work of Odlyzko and Poonen on zeros of polynomials with $0, 1$ coefficients [40].

### 5.1. Arbitrary distributions.

Given $f_0(t), f_1(t), \ldots, f_n(t)$, we now ask for the expected number of real roots of the random equation

$$a_0 f_0(t) + a_1 f_1(t) + \cdots + a_n f_n(t) = 0,$$

where we will assume that the $a_i$ have an arbitrary joint probability density function $\sigma(a)$.

Define $v(t) \in \mathbb{R}^{n+1}$ by

$$v(t) = \begin{pmatrix} f_0(t) \\ \vdots \\ f_n(t) \end{pmatrix},$$

and let

$$\gamma(t) \equiv v(t)/\|v(t)\|. \tag{12}$$

Instead of working on the sphere, let us work in $\mathbb{R}^{n+1}$ by defining $\gamma(t)_\perp$ to be the *hyperplane* through the origin perpendicular to $\gamma(t)$.

Fix $t$ and choose an orthonormal basis such that $e_0 = \gamma(t)$ and $e_1 = \gamma'(t)/||\gamma'(t)||$. As we change $t$ to $t + dt$, the volume swept out by the hyperplanes will form an infinitesimal wedge. (See Figure 6.)

This wedge is the Cartesian product of a two-dimensional wedge in the plane span($e_0$,$e_1$) with $\mathbb{R}^{n-1}$, the entire span of the remaining $n-1$ basis directions. The volume of the wedge is

$$||\gamma'(t)||\ dt \int_{\mathbb{R}^n \equiv \{e_0\cdot a = 0\}} |e_1 \cdot a|\sigma(a)da^n\ ,$$

where the domain of integration is the $n$-dimensional space perpendicular to $e_0$ and $a^n$ denotes $n$-dimensional Lebesgue measure in that space. Intuitively $\|\gamma'(t)\|dt$ is the rate at which the wedge is being swept out. The width of the wedge is infinitesimally proportional to $|e_1 \cdot a|$, where $a$ is in this perpendicular hyperspace. The factor $\sigma(a)$ scales the volume in accordance with our chosen probability measure.

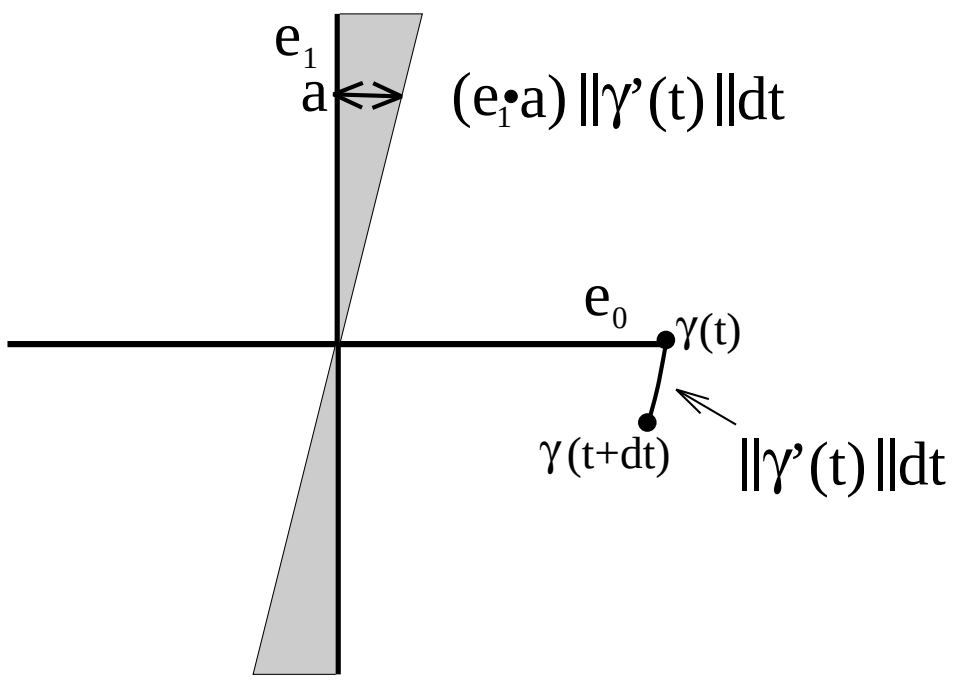

Figure 6. Infinitesimal wedge area.

**Theorem 5.1.** *If $a$ has a joint probability density $\sigma(a)$, then the density of the real zeros of $a_0f_0(t) + \cdots + a_nf_n(t) = 0$ is*

$$\rho(t)\ =\ ||\gamma'(t)|| \int_{\gamma(t)\cdot a=0} \frac{|\gamma'(t)\cdot a|}{\|\gamma'(t)\|}\ \sigma(a)\, da^n\ =\ \int_{\gamma(t)\cdot a=0} |\gamma'(t)\cdot a|\ \sigma(a)\ da^n\ ,$$

*where $da^n$ is standard Lesbesgue measure in the subspace perpendicular to $\gamma(t)$.*

5.2. **Non-central multivariate normals: Theory.** We apply the results in the previous subsection to the case of multivariate normal distributions. We begin by assuming that our distribution has mean $m$ and covariance matrix $I$. We then show that the restriction on the covariance matrix is readily removed. Thus we assume that

$$\sigma(a) = (2\pi)^{-(n+1)/2}e^{-\sum(a-m_i)^2/2},\ \text{and}\ m = (m_0, \ldots, m_n)^T.$$

**Theorem 5.2.** *Assume that $(a_0, \ldots, a_n)^T$ has the multivariate normal distribution with mean $m$ and covariance matrix $I$. Let $\gamma(t)$ be defined as in (12). Let $m_0(t)$ and $m_1(t)$ be the components of $m$ in the $\gamma(t)$ and $\gamma'(t)$ directions, respectively. The density of the real zeros of the equation $\sum a_if_i(t) = 0$ is*

$$\rho_n(t) = \frac{1}{\pi}\|\gamma'(t)\|e^{-\frac{1}{2}m_0(t)^2}\left\{e^{-\frac{1}{2}m_1(t)^2} + \sqrt{\frac{\pi}{2}}m_1(t)\mathrm{erf}\left[\frac{m_1(t)}{\sqrt{2}}\right]\right\} .$$

For polynomials with identically distributed normal coefficients, this formula is equivalent to [2, Section 4.3C].

*Proof.* Since we are considering the multivariate normal distribution, we may rewrite $\sigma(a)$ in coordinates $x_0, \dots, x_n$ in the directions $e_0, \dots, e_n$ respectively. Thus

$$\sigma(x_0, \dots, x_n) = (2\pi)^{-(n+1)/2}e^{-\frac{1}{2}\sum(x-m_i(t))^2},$$

where $m_i(t)$ denotes the coordinate of $m$ in the $e_i$ direction. The $n$-dimensional integral formula that appears in Theorem 5.1 reduces to

$$\frac{1}{2\pi}\int_{-\infty}^{\infty} |x_1|\ e^{-\frac{1}{2}(m_0(t)^2)}e^{-\frac{1}{2}(x_1-m_1(t))^2}\ dx_1$$

after integrating out the $n-1$ directions orthogonal to the wedge. From this, the formula in the theorem is obtained by direct integration. □

We can now generalize these formulas to allow for arbitrary covariance matrices as we did with Theorem 3.1. We phrase this corollary in a manner that is self-contained: no reference to definitions anywhere else in the paper is necessary.

**Corollary 5.1.** *Let $v(t) = (f_0(t), f_1(t), \dots, f_n(t))^T$, and let $a = (a_0, \dots, a_n)$ be a multivariate normal distribution with mean $m = (m_0, \dots, m_n)^T$ and covariance matrix $C$. Equivalently consider random functions of the form $\sum a_i f_i(t)$ with mean $\mu(t) = m_0 f_0(t) + \cdots + m_n f_n(t)$ and covariance matrix $C$. The expected number of real roots of the equation $\sum a_i f_i(t) = 0$ on the interval $[a, b]$ is*

$$\frac{1}{\pi}\int_a^b \|\gamma'(t)\|e^{-\frac{1}{2}m_0^2(t)}\left\{e^{-\frac{1}{2}m_1^2(t)} + \sqrt{\frac{\pi}{2}}m_1(t)\mathrm{erf}\left[\frac{m_1(t)}{\sqrt{2}}\right]\right\} dt\ ,$$

*where*

$$w(t) = C^{1/2}v(t), \quad \gamma(t) = \frac{w(t)}{\|w(t)\|}, \quad m_0(t) = \frac{\mu(t)}{\|w(t)\|}, \quad \textit{and} \quad m_1(t) = \frac{m_0'(t)}{\|\gamma'(t)\|}.$$

*Proof.* There is no difference between the equation $a \cdot v = 0$ and $C^{-1/2}a \cdot C^{1/2}v = 0$. The latter equation describes a random equation problem with coefficients from a multivariate normal with mean $C^{-1/2}m$ and covariance matrix $I$. Since $\mu(t)/\|w(t)\| = \gamma(t) \cdot C^{-1/2}m$ and $m_0'(t)/\|\gamma'(t)\| = \gamma'(t) \cdot C^{-1/2}m/\|\gamma'(t)\|$, the result follows immediately from Theorem 5.2. □

The reader may use this corollary to compute the expected number of roots of a random monic polynomial. In this case $m = e_n$ and $C$ is singular, but this singularity causes no trouble. We now proceed to consider more general examples.

5.3. **Non-central multivariate normals: Applications.** We explore two cases in which non-central normal distributions have particularly simple zero densities:

- Case I. $m_0(t) = m$ and $m_1(t) = 0$.
- Case II. $m_0(t) = m_1(t)$.

*Case* I . $m_0(t) = m$ and $m_1(t) = 0$. If we can arrange for $m_0 = m$ to be a constant, then $m_1(t) = 0$ and the density reduces to

$$\rho(t) = \frac{1}{\pi}\|\gamma'(t)\|e^{-\frac{1}{2}m^2} .$$

In this very special case, the density function for the mean $m$ case is just a constant factor ($e^{-\frac{1}{2}m^2}$) times the mean zero case.

This can be arranged if and only if the function $\|w(t)\|$ is in the linear space spanned by the $f_i$. The next few examples show when this is possible. In parentheses, we indicate the subsection of this paper where the reader may find the mean zero case for comparison.

**Example 1** (3.1.2)**.** A random polynomial with a simple answer, even degree**:** Let $f_i(t) = t^i, i = 0, \ldots, n$, and $C = \mathrm{diag}[\binom{n}{i}]$, so that $\|w(t)\| = (1+t^2)^{n/2}$. Choose $\mu(t) = m(1+t^2)^{n/2}$, so that $m_0(t) = m$ is a constant.

For example, if $n = 2$ and $a_0, a_1$, and $a_2$ are independent standard Gaussians, then the random polynomial

$$(a_0 + m) + a_1\sqrt{2}t + (a_2 + m)t^2$$

is expected to have

$$\sqrt{2}e^{-m^2/2}$$

real zeros. The density is

$$\rho(t) = \frac{1}{\pi}\frac{\sqrt{2}}{(1+t^2)}e^{-m^2/2}.$$

Note that as $m \to \infty$, we are looking at perturbations to the equation $t^2 + 1 = 0$ with no real zeros, so we expect the number of real zeros to converge to 0.

**Example 2** (3.2.4)**.** Trigonometric sums**:** $\mu(t) = m\sqrt{\sigma_0^2 + \cdots + \sigma_n^2}$.

**Example 3** (3.2.2)**.** Random power series**:** $\mu(t) = m(1-t^2)^{-1/2}$.

**Example 4** (3.2.3)**.** Entire functions**:** $\mu(t) = me^{t^2/2}$.

**Example 5** (3.2.5)**.** Dirichlet series**:**

$$\mu(t) = m\sqrt{\zeta(2t)} = \sum_{k=1}^{\infty}\frac{m_k}{k^t},$$

where $m_k = 0$ if $k$ is not a square, and $m_k = m\prod_i \frac{(2n_i-1)!!}{(2n_i)!!}$ if $k$ has the prime factorization $\prod_i p_i^{2n_i}$.

*Case* II . $m_0(t) = m_1(t)$. We may pick a $\mu(t)$ for which $m_0(t) = m_1(t)$ by solving the first-order ordinary differential equation $m_0(t) = m_0'(t)/\|\gamma'(t)\|$. The solution is

$$\mu(t) = m\|w(t)\| \ \ \exp\left[\int_K^t ||\gamma'(x)||\,dx\right] .$$

There is really only one integration constant since the result of shifting by $K$ can be absorbed into the $m$ factor. If the resulting $\mu(t)$ is in the linear space spanned by the $f_i$, then we choose this as our mean.

Though there is no reason to expect this, it turns out that if we make this choice of $\mu(t)$, then the density may be integrated in closed form. The expected number of zeros on the interval $[a, b]$ is

$$\int_a^b \rho(t)dt = \frac{1}{4}\mathrm{erf}^2(m_0(t)/\sqrt{2}) - \frac{1}{2\pi}\Gamma[0, m_0^2(t)]\Big|_a^b .$$

**Example 6** (3.2.2)**.** Random power series**:** Consider a power series with independent, identically distributed normal coefficients. In this case $\mu(t) = \frac{m}{1-t}$, where $m =$ (mean/standard deviation), so $m_0(t) = m\sqrt{\frac{1+t}{1-t}}$. A short calculation shows that $m_1(t) = m_0(t)$.

**Example 7** (3.2.3)**.** Entire functions**:** In this case $\mu(t) = me^{t+t^2/2}$, so $m_0(t) = me^t$.

**Example 8** (3.2.5)**.** Dirichlet series**:** This we leave as an exercise. Choose $K > 1/2$.

**Theorem 5.3.** *Consider a random polynomial of degree $n$ with coefficients that are independent and identically distributed normal random variables. Define $m \neq 0$ to be the mean divided by the standard deviation. Then as $n \to \infty$,*

$$E_n = \frac{1}{\pi}\log(n) + \frac{C_1}{2} + \frac{1}{2} - \frac{\gamma}{\pi} - \frac{2}{\pi}\log|m| + O(1/n) \ ,$$

*where $C_1 = 0.6257358072...$ is defined in Theorem 2.2 and $\gamma = 0.5772156649...$ is Euler's constant. Furthermore, the expected number of positive zeros is asymptotic to*

$$\frac{1}{2} - \frac{1}{2}\mathrm{erf}^2(|m|/\sqrt{2}) + \frac{1}{\pi}\Gamma[0, m^2].$$

*Sketch of proof.* We break up the domain of integration into four subdomains: $(-\infty, -1]$, $[-1, 0]$, $[0, 1]$, and $[1, \infty)$. Observe that the expected number of zeros on the first and second intervals are the same, as are the expected number of zeros on the third and fourth intervals. Thus we will focus on the first and third interval, doubling our final answer.

The asymptotics of the density of zeros is easy to analyze on [0,1] because it converges quickly to that of the power series (Example 6, above). Doubling this gives us the expected number of positive zeros.

On the interval $(-\infty, -1]$, one can parallel the proof of Theorem 2.2. We make the change of variables $-t = 1 + x/n$. The weight due to the non-zero mean can be shown to be $1 + O(1/n)$. Therefore, the asymptotic series for the density of the zeros is the same up to $O(1/n)$. We subtract the asymptotic series for the density of the zeros of the non-central random power series and then integrate term by term. □

The $\frac{1}{\pi}\log(n)$ term was first derived by Sambandham. Farahmand [17] has improved on his results.

## 6. Eigenvalues of random matrices

Eigenvalues of random matrices arise in a surprising number of disciplines of both pure and applied mathematics. Already three major books [19], [37], [38] on the subject exist, each specializing in different disciplines, yet these books serve as mere stepping-stones to the vast literature on the subject. The book by Mehta [37] covers developments of random matrices (mostly symmetric) that began with the work of Wigner, who modeled heavy atom energies with random matrix eigenvalues. Muirhead's book [38] focuses on applications to multivariate statistics, including eigenvalue distributions of Wishart matrices. These are equivalent to singular value distributions of rectangular matrices whose columns are iid multivariate normal. His exposition is easily read with almost no background. Girko's large book [19] translates his earlier books from Russian and includes more recent work as well.

An entire semester's interdisciplinary graduate course [12] was inadequate for reviewing the subject of eigenvalues of random matrices. Some exciting recent developments may be found in books by Voiculescu, Dykema, and Nica [50] relating Wigner's theory to free random variables and by Faraut and Koranyi [18], which extend the special functions of matrix argument described in [38] from the harmonic analysis viewpoint. Other new areas that we wish to mention quickly concern matrix models for quantum gravity [1], Lyapunov exponents [39], and combinatorial interpretations of random matrix formulas [20], [24]. By no means should the handful of papers mentioned be thought of as an exhaustive list.

Developers of numerical algorithms often use random matrices as test matrices for their software. An important lesson is that a random matrix should not be equated to the intuitive notion of a "typical" matrix or the vague concept of "any old" matrix. Random matrices, particularly large ones, have special properties of their own. Often there is little more information obtained from 1,000 random trials than from one trial [14].

6.1. **How many eigenvalues of a random matrix are real?** Assume that we have a random matrix with independent standard normal entries. If $n$ is even, the expected number of real eigenvalues is

$$E_n = \sqrt{2} \sum_{k=0}^{n/2-1} \frac{(4k-1)!!}{(4k)!!},$$

while if $n$ is odd,

$$E_n = 1 + \sqrt{2} \sum_{k=1}^{(n-1)/2} \frac{(4k-3)!!}{(4k-2)!!}.$$

As $n \to \infty$,

$$E_n \sim \sqrt{2n/\pi}.$$

This is derived in [13] using zonal polynomials. The random eigenvalues form an interesting Saturn-like picture in the complex plane. Figure 7 plots normalized eigenvalues $\lambda/\sqrt{50}$ in the complex plane for fifty matrices of size $50 \times 50$. There are 2,500 dots in the figure. Girko's [19] circular law (which we have not verified) states under general conditions that as $n \to \infty$, $\lambda/\sqrt{n}$ is uniformly distributed on the disk. If the entries are independent standard normals, a proof may be found in [15], where also may be found a derivation of the repulsion from the real axis that is clearly visible.

Girko's circular law stands in contrast to the result that roots of random polynomials are uniformly distributed on the unit circle rather than the disk.

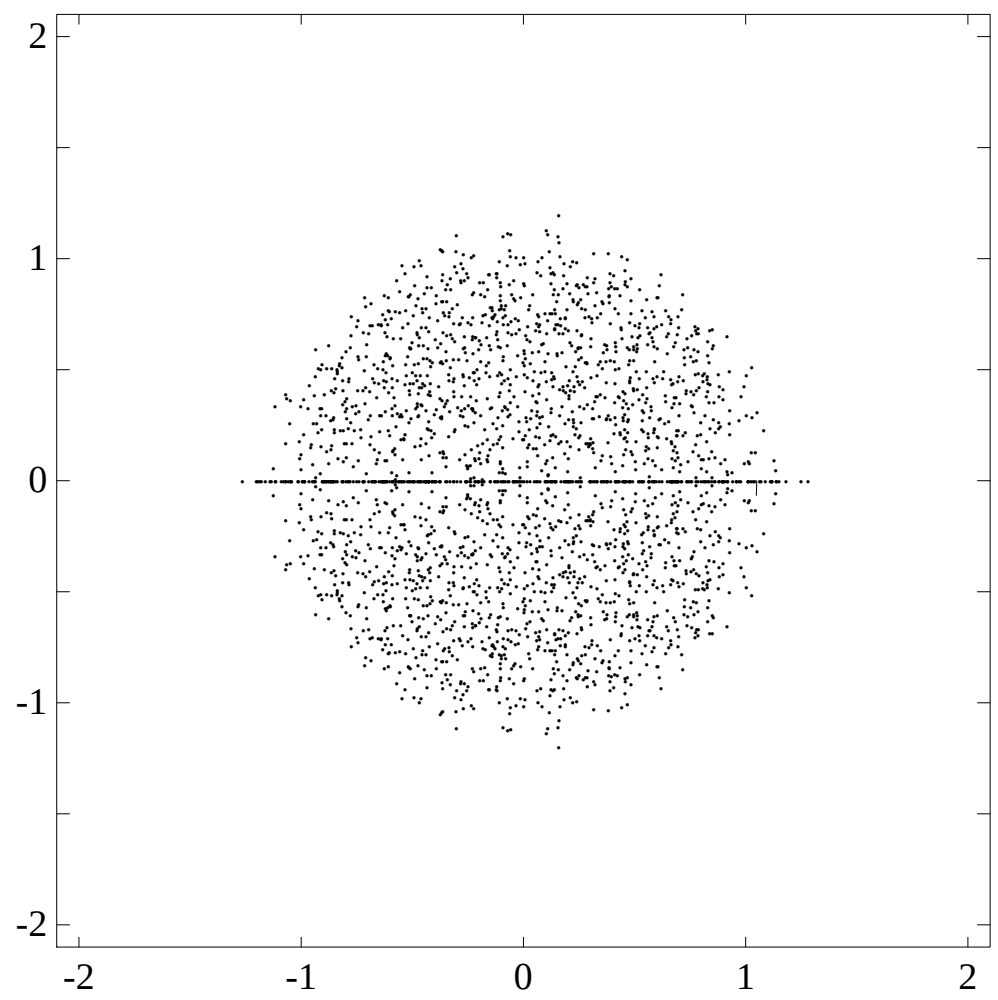

FIGURE 7. 2,500 dots representing normalized eigenvalues of fifty random matrices of size $n = 50$. Clearly visible are the points on the real axis.

6.2. **Matrix polynomials.** This may come as a shock to some readers, but characteristic polynomials are a somewhat unnatural way to discuss the eigenvalues of a matrix. It seems irrelevant that a random matrix happens to have a random characteristic polynomial, so we will not discuss random characteristic polynomials any further. An analogous situation occurs in the numerical computation of eigenvalues, where nobody would dream of forming the characteristic polynomial.

The proper generalization that includes polynomials and matrices as special cases is the so-called *matrix polynomial.* A matrix polynomial has the form

$$P(t) = A_0 + A_1 t + A_2 t^2 + \cdots + A_n t^n,$$

where the $A_i$ are $p \times p$ matrices and $t$ is a scalar. The solutions to $\det P(t) = 0$ are the eigenvalues of the matrix polynomial. Notice that we are not trying to set $P(t)$ to be the zero matrix, but rather we are trying to find a $t$ for which $P(t)$ is a singular matrix. It is sometimes convenient to take $A_n = I$. The standard eigenvalue problem takes $n = 1$ and $A_1 = I$. When $n = 1$ and $A_1 \neq I$, the problem is known as the *generalized eigenvalue problem*. Pure polynomials correspond to $p = 1$.

The beauty of random matrix polynomials is that the expected number of real eigenvalues depends on $p$ by a geometric factor:

**Theorem 6.1.** *Let $f_0(t), \ldots, f_n(t)$ be any collection of differentiable functions, and let $A_0, \ldots, A_n$ be $p \times p$ random matrices with the property that the $p^2$ random vectors $((A_0)_{ij}, (A_1)_{ij}, \ldots, (A_n)_{ij})$ $(i, j = 1, \ldots, p)$ are iid multivariate normals with mean zero and covariance matrix $C$. Let $\alpha_p$ denote the expected number of real solutions in the interval $[a, b]$ to the equation*

$$0 = \det\left[A_0 f_0(t) + A_1 f_1(t) + \cdots + A_n f_n(t)\right].$$

*We then have that*

$$\alpha_p/\alpha_1 = \sqrt{\pi}\,\frac{\Gamma((p+1)/2)}{\Gamma(p/2)}.$$

*$\alpha_1$ may be computed from Theorem 3.1.*

In particular, if all of the matrices are independent standard normals, the expected number of real solutions is

$$E_n \times \sqrt{\pi}\frac{\Gamma((p+1)/2)}{\Gamma(p/2)},$$

where $E_n$ is the quantity that appears in Theorem 2.1. The proof of Theorem 6.1 follows from a simple consequence of the integral geometry formula.

**Lemma 6.1.** *Choose an interval $[a, b]$ and a random function*

$$a_0 f_0(t) + a_1 f_1(t) + \cdots + a_n f_n(t), \qquad t \in [a, b],$$

*where the $a_i$ are independent standard normals. Generate a random curve in $\mathbb{R}^k$ by choosing an independent sample of $k$ such functions. The expected length of the projection of this curve onto the unit sphere in $\mathbb{R}^k$ is equal to the expected number of zeros of the chosen random function on the chosen interval, multiplied by $\pi$.*

*Proof.* The lemma follows from (11). Let $N$ be the random curve. Since the distribution of $QN$ is the same as that of $N$, we may take $M$ to be any fixed hyperplane, say $x_1 = 0$. The intersections of a curve with this hyperplane are exactly the zeros of the first coordinate of the curve and thus the zeros of our random function. □

Notice that the expected length does not depend on $k$. This result generalizes to random embeddings of manifolds in Euclidean space. See [33] for a discussion of these and other random varieties.

*Proof of Theorem 6.1.* We prove this theorem by using integral geometry twice to obtain expressions for the average length of the random curve $\gamma$ defined by

$$A(t) = A_0 f_0(t) + A_1 f_1(t) + \cdots + A_n f_n(t),$$

on some interval $[a, b]$, and $\gamma(t) = A(t)/\|A(t)\|_F$.

On the one hand, Lemma 6.1 states that the expected length of the projection $\gamma(t)$ is $\alpha_1 \pi$.

On the other hand, (11) may be used with $M$ chosen to be the set of singular matrices on $S^{p^2-1}$, and $N$ is the random curve $\gamma$. Thus the expected number of $t$ for which $\gamma(t)$ is singular is

$$\alpha_p = \frac{1}{\pi}\frac{|M||N|}{|S^{p^2-2}|}. \tag{13}$$

The volume of $M$ is known [13] to be

$$|M| = \frac{2\pi^{p^2/2}\Gamma((p+1)/2)}{\Gamma((p/2)\Gamma((p^2-1)/2)} .$$

The average length of $N$ is $\pi\alpha_1$. The volume of $S^{p^2-2}$ is $2\pi^{(p^2-1)/2}/\Gamma((p^2-1)/2)$. Plugging these volumes back into (13) yields the result. □

## 7. Systems of equations

The results that we have derived about random equations in one variable may be generalized to systems of $m$ equations in $m$ unknowns. What used to be a curve $v(t) : \mathbb{R} \to \mathbb{R}^{n+1}$ is now an $m$-dimensional surface $v(t) : \mathbb{R}^m \to \mathbb{R}^{n+1}$ defined in the same way. The random coefficients now form an $m \times (n+1)$ matrix $A$.

**Theorem 7.1.** *Let $v(t) = (f_0(t), \ldots, f_n(t))^T$ be any differentiable from $\mathbb{R}^m$ to $\mathbb{R}^{n+1}$, let $U$ be a measurable subset of $\mathbb{R}^m$, and let $A$ be a random $m \times (n+1)$ matrix. Assume that the rows of $A$ are iid multivariate normal random vectors with mean zero and covariance matrix $C$. The expected number of real roots of the system of equations*

$$Av(t) = 0$$

*that lie in the set $U$ is*

$$\pi^{-\frac{m+1}{2}}\Gamma\left(\frac{m+1}{2}\right)\int_U \left(\det\left[\frac{\partial^2}{\partial x_i \partial y_j}\left(\log v(x)^T C v(y)\right)\Big|_{y=x=t}\right]_{ij}\right)^{1/2} dt.$$

*Proof.* This is an application of the integral geometry formula (11). To apply this formula on the unit sphere $S^n \subset \mathbb{R}^{n+1}$, we choose a submanifold $M$ of dimension $m$ and a submanifold $N$ of dimension $n-m$.

For simplicity assume first that $C = I$. We take $M$ to be the projection of $\{v(t) : t \in U\}$ to the unit sphere. For $N$ we take the intersection of a plane of dimension $n-m+1$ with the sphere, i.e., $N = S^{n-m} \subset S^n$.

According to (11), if we intersect $M$ with a random $(n-m+1)$-dimensional plane, the expected number of intersections is

$$E(\#(M \cap QN)) = \frac{2}{|S^m||S^{n-m}|}|M||N| = 2|M|/|S^m| = \pi^{-\frac{m+1}{2}}\Gamma\left(\frac{m+1}{2}\right)|M|.$$

The Fubini-Study metric conveniently tells us that $|M|$ is the integral in the statement of the theorem.

Of course, the number of real roots of $Av(t) = 0$ is the number of intersections of $M$ with the null-space of $A$ (counting multiplicity). Since for the moment we assume that $C = I$, the random null-space of $A$ is invariant under rotations, proving that the average number of intersections is the average number of real roots.

For arbitrary $C$ the entire derivation applies by replacing $A$ with $AC^{-1/2}$. □

We now extend our previous examples to random systems of equations.

7.1. **The Kac formula.** Consider systems of polynomial equations with independent standard normal coefficients. The most straightforward generalization occurs if the components of $v$ are all the monomials $\{\prod_{k=1}^m x_k^{i_k}\}$, where for all $k$, $i_k \leq d$. In other words, the Newton polyhedron is a hypercube.

Clearly,

$$v(x)^T v(y) = \prod_{i=1}^{m} \sum_{k=0}^{d} (x_i y_i)^k,$$

from which we see that the matrix in the formula above is diagonal and that the density of the zeros on $\mathbb{R}^m$ breaks up as a product of densities on $\mathbb{R}$. Thus if $E_d^{(m)}$ represents the expected number of zeros for the system,

$$E_d^{(m)} = \pi^{-\frac{m+1}{2}} \Gamma\left(\frac{m+1}{2}\right) (\pi E_d^{(1)})^m.$$

The asymptotics of the univariate Kac formula shows that as $d \to \infty$,

$$E_d^{(m)} \sim \pi^{-\frac{m+1}{2}} \Gamma\left(\frac{m+1}{2}\right) (2 \log d)^m.$$

We suspect that the same asymptotic formula holds for a wide range of Newton polyhedra, including the usual definition of degree: $\sum_{k=1}^m i_k \leq d$ [32].

7.2. **A random polynomial with a simple answer.** Consider a system of $m$ random polynomials, each of the form

$$\sum_{i_1, \dots, i_m} a_{i_1 \dots i_m} \Pi_{k=1}^m x_k^{i_k},$$

where $\sum_{k=1}^m i_k \leq d$ and where the $a_{i_1 \dots i_m}$ are independent normals with mean zero and variances equal to multinomial coefficients:

$$\binom{d}{i_1, \dots, i_m} = \frac{d!}{(d - \sum_{k=1}^m i_k)! \prod_{k=1}^m i_k!}.$$

The multinomial theorem simplifies the computation of

$$v(x)^T C v(y) = \sum_{i_1, \dots, i_m} \binom{d}{i_1, \dots, i_m} \prod_{k=1}^m x^{i_k} y^{i_k} = (1 + x \cdot y)^d.$$

We see that the density of zeros is

$$\rho(t) = \pi^{-\frac{m+1}{2}} \Gamma\left(\frac{m+1}{2}\right) \frac{d^{m/2}}{(1 + t \cdot t)^{(m+1)/2}}.$$

In other words, the zeros are uniformly distributed on real projective space, and the expected number of zeros is $d^{m/2}$.

Shub and Smale [46] have generalized this result as follows. Consider $m$ independent equations of degrees $d_1, \dots, d_m$, each defined as in this example. Then the expected number of real zeros of the system is

$$\sqrt{\prod_{k=1}^m d_k}.$$

The result has also been generalized to underdetermined systems of equations [31]. That is to say, we may consider the expected volume of a random real projective variety. The degrees of the equations need not be the same. The key result is as follows. *The expected volume of a real projective variety is the square root of the product of the degrees of the equations defining the variety, multiplied by the volume of the real projective space of the same dimension as the variety.* For a detailed discussion of random real projective varieties, see [33].

7.3. **Random harmonic polynomials.** Consider the vector space of homogeneous polynomials of degree $d$ in $m+1$ variables that are *harmonic*, that is, the Laplacians of the polynomials are equal to zero. If $Q$ is an orthogonal $(m+1)\times(m+1)$ matrix, then the map that sends $p(x)$ to $p(Qx)$ is a linear map from our vector space to itself, i.e., we have a representation of the orthogonal group $O(m+1)$. It is a classical result in Lie group theory that there is, up to a constant, a unique normal measure on harmonic polynomials that is invariant under orthogonal rotations of the argument. It follows that this representation is irreducible.

We outline a proof by considering the invariance of $v(x)^T Cv(y)$. Assume that for any orthogonal matrix $Q$, $v(Qx)^T Cv(Qy) = v(x)^T Cv(y)$. This implies that $v(x)^T Cv(y)$ must be a polynomial in $x\cdot x$, $x\cdot y$, and $y\cdot y$. This is classical invariant theory. For proofs and discussion of such results, see [48, Vol. 5, pp. 466–486]. We thus deduce that $v(x)^T Cv(y)$ must be of the form

$$\sum_{k=0}^{[d/2]} \beta_k (x\cdot x)^k (y\cdot y)^k (x\cdot y)^{d-2k}.$$

Setting the Laplacian of this expression to zero, we see that

$$2k(m+2d-2k-1)\beta_k + (d-2k+2)(d-2k+1)\beta_{k-1} = 0$$

and therefore that

$$\frac{\beta_k}{\beta_0} = \frac{(-1)^k d!(m+2d-2k-3)!!}{2^k k!(d-2k)!(m+2d-3)!!} .$$

Thus we see that $v(x)^T Cv(y)$ is uniquely determined (up to a constant).

From this formula we can show that the expected number of roots for a system of $m$ such random harmonic polynomial equations is

$$\left(\frac{d(d+m-1)}{m}\right)^{m/2} .$$

Because of the orthogonal invariance of these random polynomials, results hold in the generality of the polynomials in Section 7.2. Thus we may consider systems of harmonic polynomials of different degrees, or we may consider underdetermined systems, and the obvious generalizations of the above result will hold. See [32] for a detailed discussion.

7.4. **Random power series.** For a power series in $m$ variables with independent standard normal coefficients, we see that the density of zeros on $\mathbb{R}^m$ breaks up as the product of $m$ densities:

$$\rho(t) = \pi^{-\frac{m+1}{2}} \Gamma\left(\frac{m+1}{2}\right) \prod_{k=1}^{m} \frac{1}{(1-t_k^2)} .$$

Notice that the power series converges with probability one on the unit hypercube, and that at the boundaries of this domain the density of zeros becomes infinite.

7.5. **Random entire functions.** Consider a random power series

$$f(x) = \sum_{i_1,\ldots,i_n} a_{i_1\ldots i_m} \Pi_{k=1}^m x_k^{i_k},$$

where the $a_{i_1\ldots i_m}$ are independent normals with mean zero and variance $\left(\prod_{k=1}^m i_k!\right)^{-1}$. Clearly

$$v(x)^T Cv(y) = e^{x\cdot y},$$

so the zeros are uniformly distributed on $\mathbb{R}^m$ with

$$\pi^{-\frac{m+1}{2}}\Gamma\left(\frac{m+1}{2}\right)$$

zeros per unit volume.

## 8. Complex zeros

We now present the complex version of Theorem 7.1 and discuss some consequences. We define a *complex* (*multivariate*) *normal vector* to be a random vector for which the real and imaginary parts are independent identically distributed (multivariate) normal vectors.

**Theorem 8.1.** *Let $v(z) = (f_0(z), \ldots, f_n(z))^T$ be any complex analytic function from an open subset of $\mathbb{C}^m$ to $\mathbb{C}^{n+1}$, let $U$ be a measurable subset of $\mathbb{C}^m$, and let $A$ be a random $m \times (n+1)$ matrix. Assume that the rows of $A$ are iid complex multivariate normal vectors with mean zero and covariance matrix $C$. The expected number of roots of the system of equations*

$$Av(z) = 0$$

*that lie in the set $U$ is*

$$\frac{m!}{\pi^m}\int_U \det\left[\frac{\partial^2}{\partial z_i \partial \bar{z}_j}\left(\log v(z)^T C v(\bar{z})\right)\right]_{ij}\prod_i dx_i\, dy_i\ .$$

*Sketch of proof.* The proof is analogous to that of Theorem 7.1 but uses complex integral geometry [43, p. 342]. The volume of the projection of $v(z)$ is calculated using the complex Fubini-Study metric [22, pp. 30–31]. □

If $U$ is Zariski open, then by Bertini's theorem, the number of intersections is constant almost everywhere. This number is called the *degree* of the embedding (or of the *complete linear system of divisors*, if we wish to emphasize the intersections). From what we have seen, the volume of the embedding is this degree multiplied by the volume of complex projective space of dimension $m$. For example, the volume of the *Veronese surface* $v : \mathbb{P}(\mathbb{C}^3) \to \mathbb{P}(\mathbb{C}^6)$, defined by

$$v(x, y, z) = (x^2, y^2, z^2, xy, xz, yz),$$

is $4 \times \pi^2/2!$. This corresponds to the fact that pairs of plane conics intersect at four points.

For the univariate case, if the coefficients are complex independent standard normals, the zeros concentrate on the unit circle (not the disk!) as the degree grows.

For the complex version of the polynomial considered in Sections 3.1.2 and 7.2, the zeros are uniformly distributed on complex projective space. Just as was observed for the real version of this example in Section 4.1, this uniformity is a consequence of (unitary) invariance of the homogeneous version of these random polynomials. But for the complex case more can be said: these polynomials provide the *unique* normal measure (up to a constant) on the space of polynomials that is unitarily invariant. A simple proof and discussion of this may be found in [30].

8.1. **Growth rates of analytic functions.** Complex analysts know that there is a connection between the asymptotic growth of analytic functions and the number of zeros inside disks of large radius. Functions whose growth may be modeled by the function $\exp(\tau z^\rho)$ are said to have *order* $\rho$ and *type* $\tau$. Precise definitions may be found in [6, p. 8]. Let $n(r)$ be the number of zeros of $f(z)$ in the disk $|z| < r$. If $f(z)$ has at least one zero anywhere on the complex plane, then [6, Eq. (2.5.19)]

$$\limsup_{r\to\infty} \frac{\log n(r)}{\log r} \le \rho. \tag{14}$$

It is possible [6, (2.2.2) and (2.2.9)] to compute the order and type from the Taylor coefficients of $f(z) = a_0 + a_1 z + a_2 z^2 + \cdots$, by using

$$\rho = \limsup_{n\to\infty} \frac{-n\log n}{\log |a_n|} \tag{15}$$

and

$$\tau = \frac{1}{e\rho} \limsup_{n\to\infty} n|a_n|^{\rho/n}. \tag{16}$$

We now illustrate these concepts with random power series. We shall restrict to the univariate case, and we shall assume that the coefficients are independent.

**Theorem 8.2.** *Let*

$$f(z) = a_0 + a_1 z + a_2 z^2 + \cdots$$

*be a random power series (or polynomial), where the $a_i$ are independent complex normals with mean zero and variances $\sigma_i^2 \ge 0$. Let*

$$\phi(z) = \sigma_0^2 + \sigma_1^2 z + \sigma_2^2 z^2 + \cdots$$

*be the generating function of the variances, and assume that $\phi(z)$ has a non-zero radius of convergence. Let $n(r)$ be the expected number of zeros of the random function $f(z)$ in the disk $|z| < r$. Then*

$$n(r) = \frac{r}{2}\frac{d}{dr} \log \phi(r^2).$$

*Proof.* Observe that $v(z)^T C v(\bar z) = \phi(z\bar z) = \phi(r^2)$, where $v(z)$ is the (infinite-dimensional) moment curve. Thus it is easy to check that

$$\frac{\partial^2}{\partial z \partial \bar z} \log v(z)^T C v(\bar z) = \frac{1}{4r}\frac{d}{dr} r \frac{d}{dr} \log \phi(r^2).$$

This is multiplied by $r dr d\theta/\pi$ and then integrated over the disk $|z| < r$. ☐

This theorem, together with the fact that the distribution of zeros is radially symmetric, completely describes the distribution of zeros for these random functions. In fact, $n(r)$ is exactly the unnormalized cumulative distribution function for the absolute values of the zeros.

As a simple example, let

$$\phi(z) = e^{2\tau z^{\rho/2}}.$$

By applying the Borel-Cantelli Lemma [45, p. 253] to (15) and (16), we see that the random function $f(z)$ has order $\rho$ and type $\tau$ with probability one. The theorem we have just established then gives

$$n(r) = \tau \rho r^\rho .$$

This result is reasonable in light of (14).

8.2. **A probabilistic Riemann hypothesis.** We conclude our discussion of complex zeros with a probabilistic analogue of the Riemann hypothesis.

**Theorem 8.3.** *Consider the random Dirichlet series*

$$f(z) = a_1 + \frac{a_2}{2^z} + \frac{a_3}{3^z} + \cdots \ , \tag{17}$$

*where $a_k$ are independent complex standard normal random variables. This converges with probability one if* $\mathrm{Re}(z) > 1/2$*. Then the expected number of zeros in the rectangle* $1/2 < x_1 < \mathrm{Re}(z) < x_2$*,* $y_1 < \mathrm{Im}(z) < y_2$*, is*

$$\frac{1}{2\pi}\left(\frac{\zeta'(2x_2)}{\zeta(2x_2)} - \frac{\zeta'(2x_1)}{\zeta(2x_1)}\right)(y_2 - y_1) \ .$$

*In particular, the density of zeros becomes infinite as we approach the critical line* $\{z \mid \mathrm{Re}(z) = 1/2\}$ *from the right.*

*Proof.* Following Section 3.2.5, we see that $v(z)^T C v(\bar{z}) = \zeta(z + \bar{z})$, so the density of zeros is

$$\frac{1}{4\pi}\frac{d^2}{dx^2}\log\zeta(2x) \ ,$$

where $x = \mathrm{Re}(z)$. □

Since (17) converges with probability one for $\mathrm{Re}(z) > 1/2$, one might try using random Dirichlet series to study the Riemann zeta function *inside* the critical strip. Unfortunately, as Section 3.2.5 and Theorem 8.3 suggest, random Dirichlet series are more closely related to $\zeta(z + \bar{z})$ than to $\zeta(z)$, and so the penetration of the critical strip is illusory.

## 9. The Buffon needle problem revisited

In 1777, Buffon showed that if you drop a needle of length $L$ on a plane containing parallel lines spaced a distance $D$ from each other, then the expected number of intersections of the needle with the lines is

$$\frac{2L}{\pi D}.$$

Buffon assumed $L = D$, but the restriction is not necessary. In fact the needle may be bent into any reasonable plane curve and the formula still holds. This is perhaps the most celebrated theorem in integral geometry and is considered by many to be the first [43].

Let us translate the Buffon needle problem to the sphere as was first done by Barbier in 1860—see [51] for a history. Consider a sphere with a fixed great circle. Draw a "needle" (a small piece of a great circle) on the sphere at random, and consider the expected number of intersections of the needle with the great circle. If we instead fix the needle and vary the great circle, it is clear that the answer would be the same.

Any rectifiable curve on the sphere can be approximated by a series of small needles. The expected number of intersections of the curve with a great circle is the sum of the expected number of intersections of each needle with a great circle. Thus the expected number of intersections of a fixed curve with a random great circle is a constant multiple of $L$, the length of the curve. To find the constant, consider the case where the fixed curve is itself a great circle. Then the average number of intersections is clearly 2 and $L$ is clearly $2\pi$. Thus the formula for the expected number of intersections of the curve with a random great circle must be

$$\frac{L}{\pi}.$$

Of course the theorem generalizes to curves on a sphere of any dimension.

To relate Barbier's result to random polynomials, we consider the curve $\gamma$ on the unit sphere in $\mathbb{R}^{n+1}$. By Barbier, the length of $\gamma$ is $\pi$ times the expected number of intersections of $\gamma$ with a random great circle. What are these intersections? Consider a polynomial $p(x) = \sum_0^n a_n x^n$, and let $\mathbf{p}_\perp$ be the equatorial $S^{n-1}$ perpendicular to the vector $\mathbf{p} \equiv (a_0, \ldots, a_n)$. Clearly $\gamma(t) \in \mathbf{p}_\perp$ for the values of $t$ where $\gamma(t) \perp \mathbf{p}$. As we saw in Section 2, these are the values of $t$ for which $p(t) = 0$. Thus the number of intersections of $\gamma$ with $\mathbf{p}_\perp$ is exactly the number of real zeros of $p$, and the expected number of real zeros is therefore the length of $\gamma$ divided by $\pi$.

## Acknowledgments

We would like to thank the referee for many valuable suggestions, Persi Diaconis for a number of interesting discussions, and John Mackey for his help in the preparation of the original draft.

Department of Mathematics Room 2-380, Massachusetts Institute of Technology, Cambridge, Massachusetts 02139

*E-mail address*: `edelman@math.mit.edu`

Arts and Sciences, Kapi'olani Community College, 4303 Diamond Head Road, Honolulu, Hawaii 96816

*E-mail address*: `kostlan@uhunix.uhcc.hawaii.edu`